\documentclass[12pt]{article}
\usepackage{amsfonts,amsmath,latexsym,amssymb,mathrsfs}
\usepackage{graphicx, url}
\usepackage{xcolor} 

\newenvironment{proof}{\noindent{{\bf Proof:}\ }}{\vspace{2ex}}

\newenvironment{prooff}{\nin{\bf Proof}}{\vspace{1ex}}

\newenvironment{remark}{{\nin\bf Remark.}}{{\vspace{2ex}}}
\def\numberlikeadb{\global\def\theequation{\thesection.\arabic{equation}}}
\numberlikeadb
\newtheorem{theorem}{Theorem}[section]
\newtheorem{lemma}[theorem]{Lemma}
\newtheorem{corollary}[theorem]{Corollary}
\newtheorem{proposition}[theorem]{Proposition}

\newtheorem{remarka}[theorem]{Remark}

\newcommand{\RR}{{\mathbb R}}

\newcommand{\proofbox}{\hspace*{\fill}\mbox{$\halmos$}}
\newcommand{\halmos}{\rule{1ex}{1.4ex}}

\def\nat{{\mathbb N}}

\def\half{{\textstyle{\frac12}}}
\def\shalf{{\scriptstyle{\frac12}}}

\def\quarter{{\textstyle{\frac14}}}

\def\Ref#1{(\ref{#1})}
\def\a{\alpha}
\def\b{\beta}
\def\s{\sigma}
\def\st{\s}
\def\f{\phi}

\def\l{\lambda}

\def\p{\pi}

\def\dtv{d_{TV}}

\def\ep{\hfill $\proofbox$ \bigskip}

\def\re{\RR}
\def\giv{\,|\,}

\def\Po{{\rm Po\,}}
\def\non{\nonumber}
\def\th{\theta}
\def\e{\varepsilon}
\def\eps{\e}

\def\var{{\rm Var\,}}
\def\r{\rho}

\def\t{\tau}
\def\g{\gamma}
\def\h{\eta}
\def\f{\phi}
\def\ps{\psi}
\def\lnti{{\lim_{n\to\infty}}}
\def\lti{{\lim_{t\to\infty}}}

\def\Bi{{\rm Bi\,}}

\def\nin{\noindent}
\def\D{\Delta}
\def\ex{{\mathbb E}}
\def\pr{{\mathbb P}}
\def\msk{\medskip}

\def\d{\delta}

\def\za{\zeta}
\def\Eq{\ =\ }
\def\Le{\ \le\ }

\def\TT{^\top}
\def\ue{r}
\def\sX{X^*}
\def\sXb{X^{(*,b)}}
\def\ur{^{(r)}}

\def\urt{^{(r,\snt)}}

\def\ui{^{(1)}}
\def\uth{^{(3)}}
\def\uf{^{(4)}}
\def\uj{^{(j)}}

\def\hY{{\widehat Y}}

\def\hPi{{\widehat \Pi}}
\def\hpi{{\hat\p}}
\def\hp{{\hat p}}

\def\tT{{\widetilde T}}
\def\tZ{{\widetilde Z}}

\def\tsi{{\tilde \s}}
\def\tz{{\tilde z}}
\def\tc{{\tilde c}}

\def\uii{^{(i)}}

\def\bals{\begin{align*}}
\def\eals{\end{align*}}
\def\bali{\begin{align}}
\def\eali{\end{align}}

\def\tQ{{\widetilde Q}}

\def\tPi{{\widetilde \Pi}}
\def\tX{{\widetilde X}}
\def\tS{{\widetilde S}}
\def\tV{{\widetilde V}}
\def\bb{\nu}
\def\bet{\beta}
\def\nuu{\zeta}
\def\be{\mathbf{e}}

\def\ignore#1{}

\def\bone{{\bf{1}}}

\def\integ{{\mathbb Z}}
\def\Z{\integ}

\def\tDe{{\widetilde \D}}
\def\hDe{{\widehat \D}}
\def\hA{{\widehat A}}
\def\tY{{\widetilde Y}}
\def\Def{\ :=\ }

\def\ut{^{(2)}}

\def\D{\Delta}

\def\hZ{{\widehat Z}}

\def\tf{{\tilde f}}

\def\tbe{{\tilde \bet}}
\def\tpi{{\tilde \p}}
\def\tp{{\tilde p}}
\def\snt{m}

\def\Prob{\pr}
\def\mathbbm{\mathbb}

\def\E{\ex}
\def\cX{{\mathcal X}}
\def\cS{{\mathcal S}}
\def\andy{\quad\mbox{and}\quad}
\def\pp{{\mathfrak p}}
\def\tpp{{\tilde{\mathfrak p}}}
\def\ujoto{^{(j_0,\snt_0)}}
\def\tYh{\tY^{(h)}}

\def\Yh{Y^{(h)}}

\def\Jh{J^{(h)}}
\def\ud{^{(d)}}
\def\ub{^{(b)}}
\def\bsb{b}

\def\textcite{\cite}
\def\hX{{\widehat X}}
\def\hq{{\hat q}}
\def\Xss#1{X_{[#1]}}

\def\cst{\g/4}
\def\cF{\mathcal{F}}
\def\cE{\mathcal{E}}

\def\hcE{\widehat{\mathcal{E}}}
\def\tal{{\tilde\a}}
\def\adbp{}
\def\adbb{}

\newcommand{\beas}{\begin{eqnarray*}}
\newcommand{\enas}{\end{eqnarray*}}
\newcommand{\eqs}{\begin{eqnarray*}}
\newcommand{\ens}{\end{eqnarray*}}

\newcommand{\bea}{\begin{eqnarray}}
\newcommand{\eqa}{\begin{eqnarray}}
\newcommand{\ena}{\end{eqnarray}}
\newcommand{\eq}{\begin{equation}}
\newcommand{\en}{\end{equation}}

\begin{document}

\title{The expected degree distribution in transient duplication divergence models}

\author{
\renewcommand{\thefootnote}{\arabic{footnote}}
A.\ D.\ Barbour\footnotemark[1]\quad
\& Tiffany Y.\ Y.\ Lo\footnotemark[2]
\\
Universit\"at Z\"urich \& University of Melbourne
}

\footnotetext[1]{Institut f\"ur Mathematik, Universit\"at Z\"urich,
Winterthurertrasse 190, CH-8057 Z\"URICH; e-mail {\tt a.d.barbour@math.uzh.ch}
}
\footnotetext[2]{School of Mathematics and Statistics, University of Melbourne, Parkville, VIC 3010, Australia;
e-mail {\tt tiffanyloyinyuan@gmail.com}
}

\date{}
\maketitle

\begin{abstract}
We study the degree distribution of a randomly chosen vertex in a duplication--divergence graph,
under a variety of different generalizations of the basic model of \cite{bhan2002duplication} and 
\cite{vazquez2003modeling}.
We pay particular attention to what happens when a non-trivial proportion of the vertices have
large degrees, establishing a central limit theorem for the logarithm of the degree distribution.
Our approach, as in \cite{jordan} and \cite{hermann2019partial}, relies heavily on the analysis of
related birth--catastrophe processes, and couplings are used to show that a number of different
formulations of the process have asymptotically similar expected degree distributions. 
\end{abstract}

\nin {\bf Keywords:}  Duplication--divergence graph, degree distribution, central limit theorem,
birth--catastrophe process.

\msk\nin
{\bf MRC subject classification:} 92C42; 05C82, 60J28, 60J85

\bals
\end{align*}

\vspace{-2cm}

\section{Introduction}\label{intro}
\setcounter{equation}{0}

   The duplication--divergence (DD) model introduced in \cite{bhan2002duplication}, \cite{vazquez2003modeling}
and \cite{chung2003duplication}, 
which we refer to as the basic model, can be described as follows. 
An initial network with~$\snt_0$ vertices is given.  At each discrete time step, a vertex in the network
is randomly chosen, and a new vertex having an identical set of neighbours is added.  With probability~$1-q$,
the new set of edges is then independently thinned, with each edge having probability~$p$ of being retained.  
Thus the new vertex has as neighbours a subset of the
neighbours of the originally chosen vertex, and hence has smaller (or equal) degree;  at the same time,
the degrees of some of its neighbours are increased by~$1$.  We henceforth refer to this model as the basic DD model.
There are other variants; for instance, a link may
be allowed between the originally chosen vertex and its copy, and the new vertex may be allowed a random
number of extra links to other, randomly chosen vertices.  Relatively little is known about the properties
of the models, though simulations suggest that they generate networks looking more like protein interaction
networks than the standard models; see, for example, \cite{vazquez2003modeling}, \cite{middendorf2005inferring} 
and \cite{gibson2011improving}.  
\adbb{
A closely related model, introduced by \cite{jordan}, is obtained by supposing that vertices 
are duplicated at random times, chosen in such a way that each vertex has an exponentially distributed `lifetime' 
before duplication, independently of all others.  The original model can then be deduced as the sequence of
states of this model at the successive duplication times.  An advantage of the model in \cite{jordan} is that other
events, such as removal of edges, can more naturally be introduced, as in \cite{hermann2019partial};  
in the original model, there is no obvious way of incorporating events that happen between duplications.
}

The information most easily deduced from the basic DD model model is the distribution of the degree of a randomly 
chosen vertex, \adbp{though other features have been addressed;  for example, the numbers of cliques and stars of different
sizes in \cite{hermann2016large}, and the evolution of the degree of a given vertex in \cite{hermann2016large} 
and \cite{frieze2020}.}
\adbb{
Starting from a configuration on~$\snt_0$ vertices, let~$N_{\snt,k}$ denote the number of vertices
of degree~$k$ in the graph when the total number of vertices is~$\snt$, and let $n_{\snt,k} := \ex N_{\snt,k}$
denote its expectation.  Then, for each $\snt \ge \snt_0$, the quantities $\pp_{\snt,k} := \snt^{-1}n_{\snt,k}$, $k \ge 0$,
define a probability distribution~$\pp_\snt$ on the non-negative integers~$\Z_+$.  Furthermore, 
$\pp_\snt$ can be interpreted as the probability distribution at time~$\snt$ of a {\it time inhomogeneous\/} Markov
chain~$Y$ on~$\Z_+$.  The state~$Y_\snt$ of the chain at time~$\snt$ represents the degree of a `tagged' vertex, changing
in response to the evolution of the DD model, with the tag switching to the {\it copy\/} whenever the tagged
vertex is chosen for duplication.  This discrete time inhomogeneous chain can be analyzed directly.  However, it was
observed in \cite{jordan} that there is a continuous time {\it homogeneous\/} pure jump Markov process~$X$ that apparently
has very similar behaviour; indeed, $X$ and~$Y$ have the same jump chain.  The process~$X$ is a birth--catastrophe
process of a particular kind, and this enables well known probabilistic techniques to be applied.
In particular, $X$ has an absorbing state at zero, with the set~$\nat$ of positive integers forming a single (transient)
communicating class, and a criterion in terms of the parameters $p$ and~$q$ can be established, to determine
whether or not the probability of absorption in zero is equal to or less than~$1$.  If the probability of absorption
in zero is~$1$, the expected proportion of vertices that have degree zero tends to~$1$ as $t \to \infty$, so
that the DD model generates a graph, most of whose vertices are isolated.  If the probability of absorption in
zero is less than~$1$, then there is an asymptotically non-trivial proportion of the vertices that have large
degrees, of magnitude growing with time.  These results, first established in \cite{hermann2016large} in the case $q=0$, 
are rather
disappointing for the practical application of the DD model, where a limiting `steady state' degree distribution would
be much preferable.
}

\adbb{
However, the degree weighted distribution, with elements 
\[
   \tpp_{\snt,k} \Def \frac{k\pp_{\snt,k}}{\sum_{l \ge 1} l\pp_{\snt,l}},
\]
can be represented as the probability distribution at time~$\snt$ of another inhomogeneous Markov chain,
to which there is also a closely related continuous time homogeneous pure jump Markov process,
studied in \cite{jordan}.  This process
evolves on~$\nat$, which is a single communicating class.  Thus, when the degree weighted process is positive recurrent, 
it has a non-trivial steady state distribution, that can be easily translated (by removing the weights)
into the distribution of the original process, conditional on non-absorption.  This, in turn, can be used as a template
for the degree distributions of networks observed in practice.  Again, a criterion for positive recurrence 
in terms of $p$ and~$q$ can be explicitly given.  What is more, a power law decay of the tail of this
distribution, with an explicit exponent, was presaged for $q=0$ in \cite{jordan}, and made precise in \cite{jacquet2020power}.
}

\adbb{
In this paper, we have two main aims.  The first is to show that the results above are relatively robust, in
that the rigid prescription of the DD model can be relaxed somewhat, without changing the qualitative
behaviour.  For instance, the thinning of edges in the copied vertex need not be binomial.  
It is relatively easy to modify the `tagged' processes, but not all such modifications have obvious
counterparts in the original DD process.  The second aim is to examine the distribution of the large degrees
when the original process is transient.  In practice, it need not be the case that an observed network is in
equilibrium, and the degree distribution of a transient process at a particular time~$\snt$ may still provide a
reasonable template.  Here, we prove a central limit theorem for the {\it logarithm\/} of the degree of a
randomly chosen vertex, if the degree is positive, implying an asymptotically log-normal distribution for the degree itself,
together with a point mass on zero.  We also make the connection between the discrete and continuous time models explicit,
showing that the discrete time model can be represented as a random non-linear time change of the continuous time model,
where the magnitude of the time change, in the time scale of the continuous process, converges almost surely as $t \to \infty$.
We conduct our analysis in the setting in which~$q$ is not necessarily equal to zero, and this also reveals 
interesting behaviour in the degree weighted process.  If $q=0$, the criterion for positive recurrence 
of the degree weighted process takes the form $p < p_o := e^{-1}$, as in \cite{jordan}.  If $0 < q < q_0 := 1/(1 + e^2)$, 
the criterion is of the
form $c_1(q) < p < c_2(q)$, where, somewhat surprisingly, $c_1(q) > 0$; above the value~$q_0$, the degree
weighted process is transient for all values of~$p$. 
}

\subsection{Model and results}\label{ADB-the-model}
 \adbb{
We work primarily with the following generalization of the DD model.  If a vertex chosen for duplication has 
degree zero, then so has its copy.
}
If a vertex chosen for duplication has degree~$k \ge 1$ and set of neighbours~$\cS$, 
we suppose that the new vertex has the same set of neighbours~$\cS$ with probability~$q_k + (1-q_k)\p_{kk}$, 
where $0 \le q_k < 1$, and that 
its neighbours consist of a random subset of~$\cS$ of size~$j$ with probability $(1-q_k)\p_{kj}$, $0 \le j < k$.
Here, for each $k \ge 1$, $\Pi_k := \bigl\{ \p_{kj},\,0\le j\le k\bigr\}$ is a probability distribution on $\{0,1,\ldots,k\}$ 
such that $\p_{kk} < 1$; \adbb{for $k=0$, $\p_{00}=1$. } 
Given that the size of the subset is~$j$, the set of neighbours is chosen uniformly among all $j$-subsets of~$\cS$.
In the basic model, $q_k = q$ for all~$k$, and~$\Pi_k$ is the binomial distribution $\Bi(k,p)$.
We then assume that the parameters are such that, at
least asymptotically as $n\to\infty$, the general model resembles the basic model.  Specifically,
we assume that, for $0 < p < 1$ and $0 \le q < 1$, and for positive constants $c_l$ and~$\g_l$, $1\le l\le 3$,
\eq\label{ADB-parameter-assns}
    |q_k - q| \Le c_1 k^{-\g_1};\quad  |p_k - p| \Le c_2 k^{-\g_2} \andy k^{-2}\s^2_k \Le c_3 k^{-\g_3},
\en 
where $kp_k$ and~$\s^2_k$ are the mean and variance of~$\Pi_k$.
\adbb{
Note that, in order to obtain a satisfactory `tagged' process, there needs to be a relationship between $p_l$ and~$q_l$,
in that  $\a_l := q_l + p_l(1-q_l)$ has to be equal to a fixed value~$\a$ for all~$l$. This is because} the probability that
the new vertex retains a given link needs to have the same value~$\a$, whatever the degree of the vertex being copied.
In the basic model, $c_1=c_2=0$, $\g_3 = 1$ and $c_3 = p(1-p)$.

Let~$N_{\snt,k}$ and $n_{\snt,k}$ be as above. 
At time $\snt+1$, a vertex~$v$ is chosen at random to be copied; with probability $\snt^{-1}N_{\snt,j}$, it has degree~$j$.
Hence the probability that the new vertex has degree~$k$, conditional on the values $(N_{\snt,l},\,l\ge0)$, is given by
\[
            \frac1{\snt}\, \Bigl\{q_k N_{\snt,k} + \sum_{j \ge k} N_{\snt,j}(1-q_j) \p_{jk} \Bigr\},
\]
and, if the new vertex has degree~$k$, this
adds one to~$N_{\snt,k}$ to be included in the count~$N_{\snt+1,k}$.  However, for any~$j$, $N_{\snt+1,j}$ may also differ from
$N_{\snt,j}$ because vertices of degree~$j$ at time~$\snt$ may become neighbours of the new vertex, each removing one from the
count of $j$-vertices at time $\snt+1$, and vertices of degree~$j-1$ may also become neighbours of the new vertex, each 
adding one to the count of $j$-vertices at time $\snt+1$.   The probability that one of the neighbours of a vertex of degree~$j$
is copied is $j/\snt$, and the probability that this event leads to a new link to the vertex is~$\a$
\adbb{(irrespective of the degree of the neighbour)}. Hence, and 
setting $n_{\snt,-1} = 0$ for all~$\snt$, it follows that
\eqa
    n_{\snt+1,k} &=& n_{\snt,k} - \frac{k\a}\snt\, n_{\snt,k} + \frac{(k-1)\a}\snt\, n_{\snt,k-1} 
  + \frac{q_k}\snt\, n_{\snt,k}
           + \frac{1}\snt\, \sum_{j \ge k} n_{\snt,j}(1-q_j) \p_{jk},\phantom{XX} \label{ADB-expected-numbers}
\ena
for  $k = 0,1,\ldots$ and $\snt \ge \snt_0$.  

\adbb{Define~$Q$ to be the matrix with elements
\begin{align}
   &Q_{k,k+1} \Eq \a k;\qquad Q_{k,k} \Eq -\{\a k + (1-q_k)(1-\p_{kk})\}; \non\\
   &Q_{k,j} \Eq  (1-q_k)\p_{kj},\quad 0\le j\le k-1;\qquad Q_{k,j} \Eq 0,\quad j \ge k+2,  \label{ADB-Q-matrix}
\end{align}
for each $k \ge 1$, and with $Q_{0,j} = 0$ for all~$j$; it is a $Q$-matrix, in the usual
Markovian sense.}
Writing $\pp_{\snt,k} := \snt^{-1}n_{\snt,k}$, Equation~\Ref{ADB-expected-numbers} gives
\eqa
    \pp_{\snt+1,k} &=&   \pp_{\snt,k} + \frac1{\snt+1}\,\Bigl\{ - \pp_{\snt,k}(1 + \a k) + \a(k-1) \pp_{\snt,k-1}   \non\\
      &&\qquad\qquad\qquad\qquad\quad\mbox{}  + q_k \pp_{\snt,k}
           +  \sum_{j \ge k} (1 - q_j)\pp_{\snt,j} \p_{jk} \Bigr\},  \label{ADB-probability-equation}
\ena
or, in vector notation,
\eq\label{ADB-Q-version-discrete}
     \pp_{\snt+1}\TT \Eq \pp_\snt\TT\{I + (\snt+1)^{-1}[Q]_{\snt+1}\}, \quad \mbox{and so}\quad 
          \pp_\snt\TT \Eq \pp_{\snt_0}\TT \prod_{j=\snt_0+1}^{\snt} \{I + j^{-1}[Q]_j\},
\en
where the matrix $[Q]_j$ consists of the rows of~$Q$ with indices $0,1,\ldots,j-1$, with all remaining rows
set to be identically zero.
The truncation of~$Q$ has no effect in~\Ref{ADB-Q-version-discrete}, because the vector~$\pp_\snt$ can only 
have positive mass on the elements $0,1,\ldots,\snt$, but it is convenient, because the
matrix $P\uj := \{I + j^{-1}[Q]_j\}$ is a stochastic matrix for all~$j$;  the row sums are all~$1$, and the 
diagonal elements are non-negative, because of the truncation.  This implies that the vector of expected
proportions~$(\pp_\snt,\,\snt\ge \snt_0)$ evolves as the state probabilities of a discrete time, inhomogeneous Markov chain, 
as in~\Ref{ADB-Q-version-discrete}: $\pp_{\snt_0+s}\TT = \pp_{\snt_0}\TT \prod_{j=1}^s P^{(\snt_0+j)}$.

Rather than work immediately with this discrete time chain, to determine the properties of~$\pp_\snt$, 
we observe that, in view of
Equation~\Ref{ADB-Q-version-discrete}, $\pp_\snt\TT$ should look much like $\pp_{\snt_0}^T e^{Q \log (\snt/\snt_0)}$,
or, alternatively, like $\pp_{\snt_0}^T P(\log \snt - \log \snt_0)$, where~$P(\cdot)$ denotes the semigroup of transition 
matrices generated by~$Q$.  Thus it is reasonable to suppose that the state probabilities at time~$\log (\snt/\snt_0)$ 
of the pure jump Markov process~$X$ 
corresponding to~$P$, with initial distribution~$\pp_{\snt_0}$, should approximate the vector~$\pp_\snt$.
\adbb{
This continuous time chain was introduced in \cite{jordan}, to describe the evolution of
the degree of a `tagged' vertex, in a model analogous to the discrete model, in which each vertex independently
has an exponential waiting time before duplication, with mean~$1$.
}
In this paper, we concentrate on the analysis of the Markov process~$X$, which, being time homogeneous, is
easier to work with than the inhomogeneous discrete time process, and derive results for the discrete process
as a consequence, in Section~\ref{ADB-discrete-time}. 
\adbb{
Since $Q_{0,0} = 0$ and, for all $k \ge 1$, $Q_{k,k+1} > 0$ and $\sum_{j=0}^{k-1}Q_{k,j} > 0$, the state zero
is absorbing for~$X$, and~$\nat$ is a single class, which is transient, because $Q_{1,0} > 0$.
Hence~$X$ is either absorbed in zero, or drifts to infinity. 
}
Letting~$\cX$ denote the event $\{X_t \to \infty\}$,
we give conditions in Section~\ref{ADB-absorption} that determine whether~$\pr[\cX]$ is equal to or greater than~$0$. 
Our arguments are based on the Foster--Lyapunov--Tweedie approach, and have the flavour of those in \cite{jordan}. 
If $\pr[\cX] = 0$, the 
{\it conditional\/} distribution of~$X_t$, given $\{X_t > 0\}$, may have a non-trivial limit, which
can also be investigated using \adbb{the degree weighted} Markov chain.   
The corresponding results, and more detailed conclusions, can be found in \cite{hermann2016large} and \cite{jordan}.
In particular, in \cite{hermann2016large}, the rate at which $\pr[X_t = 0]$ converges to its limit is investigated, also
in the case where $\pr[\cX] > 0$.
In Section~\ref{ADB-infinity}, when $\pr[\cX] > 0$, we investigate how the process~$X$ approaches infinity on the
event~$\cX$, proving an almost sure convergence theorem, reminiscient of that for supercritical branching processes.
Under these circumstances, a central limit theorem, analogous to that in \cite{bansaye2013extinction}, is established
in Section~\ref{ADB-sect-CLT}, indicating that the expected proportions of vertices of large degree in the 
corresponding DD model \adbb{follow a log-normal distribution}.

\adbb{
In Section~\ref{ADB-discrete-time}, we show that the original, discrete time model, after a logarithmic time change, 
can be closely coupled to the continuous time model. Then, in Section~\ref{ADB-immigration}, we discuss some 
other variants of the models.  In \cite{hermann2019partial}, random deletion of edges is introduced into a model
analagous to that of \cite{jordan}, in which vertices are duplicated independently of one another, each with a rate
$(\snt+1)/\snt$ if the current number of vertices is~$\snt$.  This results in transitions corresponding to deaths 
in the related (inhomogeneous) Markov process describing the tagged particle.  Here, we investigate the 
corresponding modification of the homogeneous Markov process~$X$.  We also allow multiple births to occur
simultaneously, so that, between catastrophes, the process~$X$ resembles a Markov branching process.
This may not seem natural in the DD context, though it turns out to be a useful generalization when considering
random re-wiring variants, as introduced in \cite{pastor2003evolving}.  These Markov branching catastrophe
processes are examined in their own right in \cite{hermann2020markov}.  Their analysis is based on a delightful
duality between branching catastrophe processes and a piecewise deterministic Markov process.  We use rather more
elementary techniques, which are nonetheless able to cope with the extra variation in the basic model that we allow.
There is also a large literature on birth, death and catastrophe processes, conducted under different
sets of assumptions
(see \cite{kaplan1975branching}, \cite{pakes1986}, \cite{brockwell1986} and \cite{barto1989}, for example),
none of which fit our setting; instead, we have used our own arguments.
}

\section{Absorption in zero}\label{ADB-absorption}
\setcounter{equation}{0}

In the next three sections, we study the behaviour of the pure jump Markov process
$X := (X_t,\,t\in\re_+)$ on the non-negative integers~$\Z_+$ that has $Q$-matrix
\begin{align}
   &Q_{k,k+1} \Eq k\a_k ;\qquad Q_{k,k} \Eq -\{k\a_k + \bet_k(1-\p_{kk})\}; \non\\
   &Q_{k,j} \Eq  \bet_k\p_{kj},\quad 0\le j\le k-1; \qquad Q_{k,j} \Eq 0,\quad j \ge k+2,  \label{ADB-Q-matrix-1}
\end{align}
for each $k \ge 1$, with $\a_k,\b_k > 0$, and
with $\Pi_k := \bigl\{ \p_{kj},\,0\le j\le k\bigr\}$ a probability distribution on $\{0,1,\ldots,k\}$
such that $\p_{kk} < 1$.  \adbb{For $k=0$, we define $Q_{0,j} = 0$ for all~$j$, so that the state~$0$
is absorbing.}
We assume, as before, that $kp_k$ and~$\s^2_k$ are the mean 
and variance of~$\Pi_k$, and that, as in~\Ref{ADB-parameter-assns},
\eq\label{ADB-parameter-assns-again}
   |p_k - p| \Le c_2 k^{-\g_2} \andy k^{-2}\s^2_k \Le c_3 k^{-\g_3};
\en
we also suppose that 
\eq\label{ADB-parameter-assns-gen}
    |\bet_k - \bet| \Le c_1 k^{-\g_1} \andy  |\a_k - \a| \Le c_4 k^{-\g_4},
\en 
for positive constants $\a,\bet,c_1,c_4,\g_1$ and~$\g_4$, and we define
\eq\label{ADB-gamma-def}
   \g \Def \{\min_{1\le l\le 4} \g_l\}\wedge 1.
\en 
\adbb{
The process that we study here is rather more general than the one with $Q$-matrix as in~\Ref{ADB-Q-matrix},
whose parameters are constrained to have $\a_k := 1 - (1 - p_k)\bet_k = \a > 0$ constant, where $\bet_k := 1-q_k$;
for that process,} we can then take
$\g_4 = \min\{\g_1,\g_2\}$ and $c_4 = c_1 + c_2$.
Because~$X$ is dominated by a Yule process with rate $\max_{k \ge 1}\a_k$, it is non-explosive.  

We also extend~$X$ to a bivariate Markov jump process $(X,Z) := \bigl((X_t,Z_t),\,t\in\re_+\bigr)$, whose $Q$-matrix
has elements
\begin{align}
   Q_{(k,l),(k+1,l)} &\Eq k\a_k ;\qquad Q_{(k,l),(k,l)} \Eq -\{k\a_k + \bet_k\}; \non\\
   Q_{(k,l),(j,l+1)} &\Eq  \bet_k\p_{kj},\quad 0\le j\le k,  \label{ADB-Q-matrix-2}
\end{align}
\adbb{for all $k \ge 0$, with $\bet_0 := 0$; for all other pairs $(k',l')$, $Q_{(k,l),(k',l')}:=0$.}
The $X$-component yields the same process as before, and the $Z$-component is the counting process
of downward jumps, including (for neatness) a `jump' at rate $\bet_k\p_{kk}$ at which the value of~$X$ does
not change.  Thus, in state~$(k,l)$, the jump rate for~$Z$ is $\bet_k$, indicating that~$Z$ 
is not too far from being a homogeneous Poisson process with rate~$\bet$, which is exactly the case in the basic model.
We let $\cF_t := \s((X_s,Z_s),\,0 \le s \le t)$.

The state~$0$ is absorbing for~$X$.
Because of the assumptions that \adbb{$\bet_k > 0$ and $\pi_{kk} < 1$} for all~$k \ge 1$, it follows that 
the set $\nat := \{1,2,\ldots\}$
forms a single class for~$X$, which is not closed, because it is possible to reach state~$0$ from state~$1$, and hence
is transient.  As before, let~$\cX$ denote the event $\{X_t \to \infty\}$.  
The main result of this section is a criterion for determining whether $\pr[\cX]$ is equal to or greater than zero.
\adbb{
The arguments that we use, as in \cite{jordan}, are based on Foster--Lyapunov--Tweedie theorems for irreducible Markov 
processes, which are sufficiently robust to accommodate our more general model.  However, since~$X$ is not
irreducible, we first modify it to an irreducible process~$\sX$, by adding 
a transition from state~$0$ to state~$1$  with positive rate~$\a_0$.
All states in~$\Z_+$ are now either positive recurrent, null recurrent or transient for~$\sX$, since it is irreducible; 
$\Prob[\cX ]= 0$ if~$\sX$ is recurrent, and  $\Prob[\mathcal{X}] > 0$ if~$\sX$ is transient. 
A criterion for deciding which is the case is established in the following theorem.
}

\begin{theorem}\label{ADB-hX-recurrence}
For the process~$\sX$, we have the following behaviour:
 \begin{align}
 &\mathrm{(1)}:\quad \mbox{if}\ \alpha < \bet \log(1/p),\ \mbox{then~$\sX$ is geometrically ergodic};\non\\
 &\mathrm{(2)}:\quad \mbox{if}\ \alpha = \bet \log(1/p),\ \mbox{then~$\sX$ is null recurrent};\label{ADB-FLT-hat-X}\\
 &\mathrm{(3)}:\quad \mbox{if}\ \alpha > \bet \log(1/p),\ \mbox{then~$\sX$ is transient}.\non
\end{align}
\end{theorem}

\begin{proof}
To prove~(1), we apply \cite[Theorem~4.2]{meyn1993stability} with $f(j) := (j+1)^{\h}$ for $\h \in (0,1]$ suitably chosen.
\adbb{
It is enough then to show that $(Qf)(j) \le -c f(j)$ for all~$j$ sufficiently large, for some $c > 0$,
since, for all $j \ge 0$, $\sum_{l \ne j} Q_{j,l}(l+1)^{\h} < \infty$.
} 
Letting $Y_j\sim \Pi_j$,  we have
\begin{align}
   (Qf)(j) &= j\alpha_j  [f(j+1)-f(j)] 
         + \bet_j(\ex\{f(Y_j)\}-f(j)), &j \ge 1. \label{ADB-Qf-def}
\end{align}
Now
\[
     j[f(j+1)-f(j)] \Eq \h(j+1)^\h (1 + O((j+1)^{-1}),
\]
and~$f$ is concave, so that
\eq\label{ADB-eta-concave}
   \ex\{f(Y_j)\}-f(j) \Le f(\ex Y_j) - f(j) \Eq (jp_j+1)^\h - (j+1)^\h
            \Eq (j+1)^\h(p^\h-1 + O(j^{-\g})),
\en
where~$\g$ is as defined in~\Ref{ADB-gamma-def}.  Hence it follows that
\[
    (Qf)(j) \Le f(j)\{ \a\h + \bet(p^\h-1) + O(j^{-\g})\}.
\]
Now 
\adbb{
\[
      \frac d{d\h} \{\a\h + \bet(p^\h-1)\} \Eq \a - \bet p^{\h}\log(1/p) \ <\ 0
\]
for all $\h > 0$ small enough, if $\a < \bet\log(1/p)$, so that there are values of $\h \in (0,1]$
such that $\a\h < \bet(1-p^\h)$.  Then, for any such~$\h$,
\[
    (Qf)(j) \Le -c_\h f(j)
\]
for all~$j$ large enough, where $c_\h := \half\{\bet(1 - p^\h) - \a\h\} > 0$.
}
\ignore{
$e^{-x} \le 1 - x + \half x^2$ in $x > 0$, so that
\[
   \bet(1 - p^\h) \ \ge\ \bet\{\h\log(1/p) - \half(\h\log(1/p))^2\}.
\]
Thus we have
\[
    (Qf)(j) \Le -f(j)\{\h(\bet\log(1/p) - \a) - \half\bet(\h\log(1/p))^2 + O(j^{-\g})\},
\]
and, for~$\h$ chosen small enough and for $\a < \bet\log(1/p)$, 
\[
   c \Def \half\{\h(\bet\log(1/p) - \a) - \half\bet(\h\log(1/p))^2\}\ >\ 0.
\]
For this choice of~$c$,
we thus have $(Qf)(j) < -cf(j)$ for all~$j$ large enough.
} 
This proves~(1).

To prove that~$\sX$ is recurrent when $\alpha = \bet \log(1/p)$, 
we apply \cite[Theorem~3.3]{tweedie1975sufficient} to the jump chain of~$\sX$. We take 
$f(j):=\log(j+1) + 1-(j+1)^{-\nuu}$ with $0 < \nuu < \g$, where~$\g$ is as defined in~\Ref{ADB-gamma-def};
then $\lim_{j\to\infty}f(j) = \infty$, and  we show that 
${Qf}(j) \leq 0$ for all $j$ sufficiently large, if $\alpha = \bet \log(1/p)$. 
It is immediate that
\eq\label{ADB-j-diff-bdd}
             j[f(j+1) - f(j)] \Eq (1 + \nuu (j+1)^{-\nuu}) + O((j+1)^{-1})
\en
is bounded for all~$j$, and that
\[
       f'(x) \Eq \frac1{x+1} + \frac{\nuu}{(x+1)^{1+\nuu}} 
\]
is decreasing in $x > 0$, so that~$f$ is concave. Hence, applying Jensen's inequality, it follows that
\ignore{
\[
               g(x) \Def \tf(x)\bone_{[0,x_\nuu)}(x) + f(x)\bone_{[x_\nuu,\infty)}(x),
\]
and because $\var(Y_j) \le c_2 j^{2-\g}$, it follows that, for~$j$ such that $jp_j \ge x_\nuu$,
}
\eqs
      \ex\{f(Y_j)\} - f(j) &\le& f(\ex Y_j) - f(j) \Eq  \log p + j^{-\nuu}(1 - p^{-\nuu}) + O(j^{-\g}) . 
\ens
Hence, from~\Ref{ADB-Qf-def},
\eqa
   (Qf)(j) &=& j\a_j[f(j+1)-f(j)] + \bet_j\bigl\{\ex\{f( Y_j)\} - f(j)\bigr\} \label{ADB-Qf-1}\\
       &\le& j\a[f(j+1)-f(j)] + \bet\bigl\{f(\ex Y_j) - f(j)\bigr\}  + O(j^{-\g})\label{ADB-Qf-2}\\
           &=& \a (1 + \nuu j^{-\nuu}) 
             + \bet\{\log p + j^{-\nuu}(1 - p^{-\nuu})\}  + O(j^{-\g}) \label{ADB-Qf-3}\\
       &=&  j^{-\nuu}\{\a\nuu + \bet(1 - p^{-\nuu})\} + O(j^{-\g}), \non
\ena
because $\a - \bet\log(1/p) = 0$.
Since $e^x > 1+x$ in $x > 0$, it follows that $p^{-\nuu} > 1 + \nuu\log(1/p)$, and hence that
\[
        \a\nuu + \bet(1 - p^{-\nuu}) \ <\ \nuu(\a - \bet\log(1/p)) \Eq 0;
\]
thus, because $\nuu < \g$, $(Qf)(j) \le 0$ for all~$j$ large enough.

To show that~$\sX$ is null, if $\alpha = \bet \log(1/p)$, we apply 
a theorem of \textcite[Appendix]{barbour1996thresholds}, which is just the continuous time analogue of 
\cite[Theorem~9.1(ii)]{tweedie1976criteria}.    
Now taking $f(j):=\log(j+1) + (j+1)^{-\nuu}$, with $0 < \nuu < \g$ as before,
we show that $(Qf)(j) \ge 0$ for all~$j$ large enough, and also that $\sum_{k\ge0}Q_{jk}|f(j) - f(k)|$ is uniformly
bounded for all~$k$.  As above, the key step is to approximate the expectation $\ex\{f(Y_j)\}$, now with a
{\it lower\/} bound.  For this we use the lower bound
\eq\label{ADB-Jensen-lower}
     f(y) \ \ge\  f(x) + (y-x)f'(x) + (y-x)^2 f''(x),
\en
vaild for~$y$ in some neighbourhood of~$x$, in view of Taylor's theorem, and because,
since $0 < \nuu < 1$, $f''(x) < 0$ for all $x > 0$.  Indeed, because $f$, $f'$ and~$f''$ are uniformly close to~$\log x$
and its derivatives for large~$x$, it is easy to check that
\Ref{ADB-Jensen-lower} is valid for all $y \ge x/2$, provided that $x \ge x_\nuu$, for some fixed $x_\nuu > 0$.
For~$j$ such that $jp_j > x_\nuu$, we thus have
\eqa\label{ADB-Ef-lower-1}
    \ex\{f(Y_j)\} &\ge& \ex\Bigl\{ \bigl(f(\ex Y_j)  + (Y_j - \ex Y_j)f'(\ex Y_j) + (Y_j - \ex Y_j)^2 f''(\ex Y_j)\bigr)
             \bone_{[\shalf\ex Y_j,\infty)}(Y_j) \Bigr\} \non\\
         &\ge&  f(\ex Y_j)(1 - \pr[Y_j < \half \ex Y_j]) + f''(\ex Y_j)\var(Y_j).    
\ena
From~\Ref{ADB-parameter-assns} and because of the definition of the function~$f$, this immediately implies that, 
for all~$j$ large enough,
\eq\label{ADB-Ef-lower-2}
    \ex\{f(Y_j)\} \ \ge\ f(\ex Y_j) + O(j^{-\g}\log j),
\en
where $\pr[Y_j < \half jp_j]$ has been bounded using Chebyshev's inequality.
Now, arguing exactly as for \Ref{ADB-Qf-1} and~\Ref{ADB-Qf-2}, but with the current definition of~$f$, it follows
that, for all~$j$ large enough,
\[
    (Qf)(j)  \ \ge\  -j^{-\nuu}(\a\nuu + \bet(1 - p^{-\nuu})) + O(j^{-\g}\log j),
\]
and hence that~$(Qf)(j) \ge 0$ for all~$j$ large enough.  It is easier to show that $\sum_{k\ge0}Q_{jk}|f(j) - f(k)|$ 
is uniformly bounded by $2(\a + \bet\log(1/p))$, for all~$j$ large enough; the necessary estimates are essentially 
already in~\Ref{ADB-Qf-3}, and the extension to all~$j$ is immediate.  This completes the proof of~(2).

For part~(3), we use \cite[Theorem 11.3(i)]{tweedie1976criteria} applied to the jump chain of~$\sX$,
with the function $f(j) := (j+1)^{-\h}$, for suitably 
chosen~$\h > 0$; since~$f$ is a decreasing function, we simply need to show that $(Qf)(j) \le 0$ for all~$j$ sufficiently
large.  Now the key step is to upper bound the expectation $\ex\{f(Y_j)\}$ using
\eq\label{ADB-Jensen-upper}
    f(y) \Le  f(x) + (y-x)f'(x) + (y-x)^2 f''(x),
\en
since $f''(x) > 0$ for all $x > 0$, valid as for~\Ref{ADB-Jensen-lower} for all~$y \ge \th x$, where $\th = \th(\h) < 1$.
Indeed, by substituting $y = \th x$ into~\Ref{ADB-Jensen-upper}, with $f(j) = (j+1)^{-\h}$, it can be seen 
that~\Ref{ADB-Jensen-upper} is satisfied for all $y \ge \th x$, provided that~$x$ is large enough,
if~$\th$ is chosen close to~$1$ in such a way that
\eq\label{ADB-Jensen-cond}
     \th^{-\h} \ <\ 1 + \h(1-\th) + \h(1+\h)(1-\th)^2.
\en
Since the left hand side 
of~\Ref{ADB-Jensen-cond} can be written as $1 + \h(1-\th) + \half\h(1+\h)(1-\th)^2 + O(|1-\th|^3)$, there are
such values of~$\th$ to be chosen.  Now, arguing much as for \Ref{ADB-Ef-lower-1} and~\Ref{ADB-Ef-lower-2}, with
the current definition of~$f$ and with $\half$ replaced by~$\th$, it follows that
\eq\label{ADB-f-convex}
       \ex\{f(Y_j)\}  - f(j) \Le f(\ex Y_j) + O(j^{-\g})  - f(j) \Eq j^{-\h}(p^{-\h}-1) + O(j^{-\g}),
\en
with~$\g$ as before. On the other hand,
\[
     j\a_j[f(j+1) - f(j)] \Eq -\a\h j^{-\h}(1 + O(j^{-\g})),
\]
and hence, for all~$j$ large enough,
\eqs
    (Qf)(j) &=& \bigl\{-\a\h  + \bet(p^{-\h}-1)\bigr\}j^{-\h} + O(j^{-\g}) \\
          &=& \bigl\{-\a\h  + \bet(e^{\h\log(1/p)} - 1)\bigr\}j^{-\h} + O(j^{-\g}).
\ens
Thus, if $\a > \bet\log(1/p)$, it is possible to choose~$\h > 0$ small enough that
$\h < \g$ and that
\[
     -\a\h  + \bet(e^{\h\log(1/p)} - 1) \ <\ 0,
\]
in which case $(Qf)(j) \le 0$ for all~$j$ sufficiently large, establishing~(3).
This completes the proof of the theorem. \ep
\end{proof}

\begin{remarka}\label{ADB-no-zero-state}
{\rm
\adbb{
Note that, if $\p_{11} = 1$ and $\p_{k0} = 0$ for all $k \ge 1$, then the process~$X$ is irreducible on~$\nat$.
For such a process, we can take $\sX = X$, and the conclusions of Theorem~\ref{ADB-hX-recurrence}
remain valid, either using the same proof, or by re-labelling the states by subtracting~$1$ from each label.  
The degree weighted process of Section~\ref{ADB-degree-weighted} comes into this category.
}
}
\end{remarka}

As a consequence of Theorem~\ref{ADB-hX-recurrence}, we have the following criterion determining
whether $\pr[\cX]$ is zero or positive, valid because
$\Prob[\cX ]= 0$ if~$\sX$ is recurrent, and  $\Prob[\mathcal{X}] > 0$ if~$\sX$ is transient.

\begin{corollary}\label{ADB-recurrence}
 If $\a \le \bet\log(1/p)$, then $\pr[\cX]=0$;  if $\a > \bet\log(1/p)$, then $\pr[\cX] > 0$.
\end{corollary}

\begin{remark}
In the basic model, with~$\bet = (1-q)$,
the condition $\a \le \bet\log(1/p)$ is equivalent to 
\eq\label{ADB-p(q)-equation}
          pe^p \Le e^{-q/(1-q)},
\en 
and is thus satisfied for $p \le p_*(q)$, where $p_*(q)$ is the unique solution of~\Ref{ADB-p(q)-equation}
in $[0,1]$.  As in \cite{hermann2016large} and \cite{jordan}, $p_*(0) \approx 0.5671$
is the solution to $pe^p = 1$;  $p_*$ is a decreasing function, and $p_*(1) = 0$. 
\end{remark}

\adbb{
In \cite{jordan}, conditions are given as to which moments of the degree weighted process are finite.
Here, we use similar arguments to make comparable statements about the moments of~$\sX$, when it is
positive recurrent. For any $u > 1$, let $x^*(u)$ denote
the positive solution to the equation $x = u(1-e^{-x})$; note that $x^*(u) < u$.
}

\begin{theorem}\label{ADB-recurrence-2}
\adbb{
 If $\a < \bet\log(1/p)$, the process~$\sX$ has a stationary distribution, whose $\h$-th moment is finite
for all~$\h < \h^* := x^*(\bet\log(1/p)/\a)/\log(1/p)$, and is infinite for all~$\h > \h^*$.
In particular, $\h^* < \b/\a$.
}
\end{theorem}

\begin{proof}
\adbb{
If $\a < \bet\log(1/p)$, 
 $\h \le 1$ and $\a\h < \bet(1-p^\h)$, the existence of an $\h$-th moment of the stationary distribution 
of~$\sX$ follows from the proof of case~(1) in Theorem~\ref{ADB-recurrence}, by  \cite[Theorem~4.2]{meyn1993stability}.
If $\h > 1$ and $\a\h < \bet(1-p^\h)$, the inequality~\Ref{ADB-eta-concave} cannot be established as before, because
$f(j) := (j+1)^\h$ is no longer concave.  However, as for~\Ref{ADB-f-convex}, we can use the bound~\Ref{ADB-Jensen-upper},
valid for this choice of~$f$ for $|y-x| \le c(\h) x$ and for all~$x$
large enough, for some $c(\h) > 0$.  Using Chebyshev's inequality to bound any contribution from 
$|Y_j - jp_j| > c(\h)jp_j$, taking expectations in~\Ref{ADB-Jensen-upper} results in the bound 
\[
   \ex\{f(Y_j)\}-f(j) \Le f(\ex Y_j) - f(j) + O(j^{\h-\g}) 
            \Eq (j+1)^\h(p^\h-1 + O(j^{-\g})),    
\]
which is the final inequality in~\Ref{ADB-eta-concave} once again. The remainder of the argument is as before.
Note that, if $\h < \h^*$, then $\a\h < \bet(1-p^\h)$.
}

\adbb{
For $\h > \h^*$, so that $\a\h > \bet(1-p^\h)$, entirely similar arguments can be used to show that
\[
    (Qf)(j) \Eq f(j)\{ \a\h + \bet(p^\h-1) + O(j^{-\g})\},
\]
and hence that
\[
      (Qf)(j) \ \ge\ c'_\h f(j), \quad \mbox{where}\ c'_\h \Def \half\{ \a\h - \bet(1-p^\h) \} \ >\ 0,
\]
for all~$j \ge j_0$, for some~$j_0$ large enough.  Defining $f_0(j) := \max\{f(j_0),f(j)\}$, it thus
follows that
\[
      (Qf_0)(j) \ \ge\ c'_\h f_0(j), \quad j \ge j_0;\qquad (Qf_0)(j) \ \ge\ 0, \quad \mbox{otherwise}.
\]
Now, with $f(j) := (j+1)^\h$, the conditions of \cite[Theorem~2]{hamza1995conditions} are
satisfied by~$f_0$, implying that $f_0(\sX_t) - \int_{t_0}^t (Qf_0)(\sX_s)\,ds$ is a martingale in $t \ge t_0$, for
any~$t_0$.  If $\ex f_0(\sX_{t_0}) < \infty$, this in turn implies that $\ex f_0(\sX_t)$ is non-decreasing
in $t \ge t_0$, and that $\ex f_0(\sX_t) > \ex f_0(\sX_{t_0})$
for $t > t_0$ if, with positive probability, $\sX_s \ge j_0$ on a subset of $[t_0,t]$ of positive Lebesgue measure;
and this is the case if $\pr[\sX_{t_0} \ge j_0] > 0$, because $\sX$ is a pure jump Markov process.
Now the stationary distribution of~$\sX$ assigns positive probability to all states, so that this
condition is satisfied if~$\sX_{t_0}$ has the stationary distribution, in which case~$\sX_t$ also has
the stationary distribution, and $\ex f_0(\sX_t) = \ex f_0(\sX_{t_0})$ if either one is finite.
Thus $\ex f_0(\sX_{t_0}) < \infty$ is not possible, and the $\h$-th moment of the stationary distribution of~$\sX$
is infinite.
}
\ep
\end{proof}

\section{The path to infinity; almost sure convergence}\label{ADB-infinity}
\setcounter{equation}{0}

In view of Theorem~\ref{ADB-recurrence},
the event $\cX$, on which $X_t \to \infty$, has positive probability if $\a > \bet\log(1/p)$.   
We now investigate how~$X_t$ approaches infinity, when it does so.  
If~$\b_k$ took the value zero for all~$k$, and if $\a_k = \a$ were the same for all~$k$, then the
process~$X$ would be a Yule process with {\it per capita\/} birth rate~$\a$. In this case, $e^{-\a t}X_t$
would converge almost surely to a positive limit~$W$.  The presence of catastrophes, when $\b > 0$, changes the
behaviour of~$X$ substantially.  Nonetheless, there is an analogous almost sure convergence theorem for the bivariate
process $(X,Z)$, \adbb{which, for constant $\a$ and~$\bet$ and binomial thinning, would be a consequence
of \cite[Theorem~4.1]{kaplan1975branching} and \cite[Theorem~4.1]{buhler1989linear}.  The theorem} 
shows that, on~$\{\lti X_t = \infty\}$, $X_t$ 
behaves like $e^{\a t}p^{Z'_t}W$, where~$W$ is a positive random variable and~$Z'$ is a Poisson process
of rate $\bet$.  Throughout the section, we assume that 
\eq\label{ADB-beta-def}
   \bb \Def \a - \bet\log(1/p) \ >\ 0;
\en
for the basic process, $\bb = \a - (1-q)\log(1/p)$.

\begin{theorem}\label{aslimit}
If $\bb  > 0$, then $W_t := e^{-\alpha t}p^{-Z_t} X_t$ has an almost sure limit 
as $t\to\infty$, which we denote by~$W$. Furthermore, $\cX = \{W > 0\}$ almost surely.
\end{theorem}

\adbb{The process~$(W_t,\,t\ge0)$ is the counterpart, in the process with catastrophes, of the usual non-negative 
martingale appearing in the Markov branching process, though, in our setting, it is only approximately a
martingale.  Just as for branching processes, }
the random variable~$W$ is essentially determined  by the early fluctuations of the process.
It is easier to understand the behaviour by writing the convergence in the form
\eq\label{ADB-approx-growth}
       \log X_{t} - \a t \ \approx\ -Z_t\log(1/p) + \log W_s,
\en
for $s$ and~$t$ large, but $s \ll t$.  This is an approximation only because $\log W_s$ is not exactly equal to~$\log W_t$,
but the two are close, for $s$ and~$t$ large, if $W > 0$.   Thus the process $\log X_{t} - \a t$, instead of being 
asymptotically constant on $\{X_t \to \infty\}$, as in the case of the Yule process,
has downward jumps of magnitude $\log(1/p)$, approximately according to a Poisson process of rate~$\bet$, staying
very close to the lattice $\log W + \log(1/p)\integ$.  The origin of the lattice is essentially determined by~$W_s$,
and the development after~$s$ is approximately given by  
\eq\label{ADB-heuristic-2}
     (\log X_t - \a t) - (\log X_s - \a s) \ \approx\ -(Z_t - Z_s)\log(1/p),
\en
where $Z_t - Z_s$ is close to a Poisson process of rate $\bet$ that is independent of~$X_s$.

\adbb{
The proofs of the counterparts of Theorem~\ref{aslimit} in \cite[Theorem~4.1]{kaplan1975branching} and 
\cite[Theorem~4.1]{buhler1989linear} use the fact that, under their
assumptions, there is a Galton--Watson process with random environments embedded in the process, enabling techniques from
branching processes to be invoked.  In our setting, this is no longer the case.  Instead, we proceed in
a succession of intermediate steps, establishing that~$\log X_t$ develops in almost deterministically linear
fashion between catastrophes, provided that~$X_t$ is large enough.  This is already close to~\Ref{ADB-approx-growth}.
The remainder of the argument consists of tidying up the details.
}

We begin by considering the bivariate process~$(X,Z)$ starting a time~$t_0$ in an arbitrary state~$(j_0,k_0)$. \adbb{Define
\eq\label{ADB-T-defs}
    T_0 \Def t_0;\qquad T_n \Def \inf\{t > t_0\colon\, Z_t = k_0+n\},
\en 
so that~$T_n$ is the time at which~$Z$ makes its $n$-th jump after~$t_0$.}
Define $\Xss{i} := X_{T_i}$.  
With~$\g$ as in~\Ref{ADB-gamma-def}, and for $i\geq 1$, define the following events:
\begin{align}
    & A_i \Def A_i(j_0,t_0) \Def \Bigl\{ \sup_{T_{i-1}\leq t < T_i}
                      \Bigl| \log X_{t} -  \a(t-T_{i-1}) - \log X_{T_{i-1}} \Bigr| 
                           \leq X_{T_{i-1}}^{-\cst} \Bigr\}; \non\\
    & B_i \Def B_i(j_0,t_0) \Def \{ |\log X_{T_i} - \log p - \log X_{T_i-} | \leq X_{T_i-}^{-\cst} \}.
    \label{ADB-AB-defs}
\end{align}
\adbb{We begin by showing that the events $A_i$ and~$B_i$ typically hold, if~$\Xss{i-1}$ is large, as is eventually
the case if $X_t \to \infty$.}

\begin{lemma}\label{ADB-A-and-B-probs}
With~$\g$ as in~\Ref{ADB-gamma-def}, and for $i\geq 1$,
the conditional probabilities of $A_i$ and~$B_i$ can be bounded as follows:
\bals
   &(1):\quad \pr[A_i(j_0,t_0) \giv \cF_{T_{i-1}}] \Le  C_A \Xss{i-1}^{\g/2-1} \Le C_A \Xss{i-1}^{-\g/2};  \\
   &(2):\quad \pr[B_i(j_0,t_0) \giv \s(T_i \vee \cF_{T_{i}-})] \Le C_B \Xss{i-1}^{-\g/2}.
\end{align*}
Here, $C_A$ and~$C_B$ are (finite) constants depending only the sequences $(\a_j,p_j,\s^2_j,\,j\ge1)$.
\end{lemma}

\begin{proof}
From time~$T_{i-1}$ to time~$T_{i}$, $X$ evolves as a pure birth process, with jump rate~$j\a_j$ when in
state~$j$.  Define~$Y\uii_t = X_t$ for $T_{i-1} \le t < T_i$, set $Y\uii_{T_i} := Y\uii_{T_i-}$,
and continue $Y\uii_t$ for $t > T_i$ as a pure birth process with the same distribution, independently of the
further evolution of~$X$. With~$S_0 := T_{i-1}$, let~$S_l$ be the time of the $l$-th jump of~$Y\uii$, and let
\eq\label{ADB-harmonic-numbers}
      h(0) \Def 0;\qquad h(j) \Def \sum_{l=1}^j \frac1l, \qquad j\ge1,
\en
denote the harmonic numbers.  Then,  conditional on~$\cF_{T_{i-1}} \cap \{\Xss{i-1} = j\}$,  it follows that
\eq\label{ADB-M-def}
    M_l \Def h(j + l-1) - h(j-1) - \sum_{k=1}^l \a_{j+ k-1} (S_k-S_{k-1}),\qquad l \in \nat_0, 
\en
can be written as the sum $\sum_{k=1}^l E'_k$, where 
\[
      E'_k \Def (j + k-1)^{-1}(1 -E_k),
\]
and $(E_k,\,k\in\nat)$
are independent standard exponential random variables.  Note also that~$M$ can be written in the form
\eq\label{ADB-M-def-2}
            M_l \Eq h(Y\uii_{S_l}-1) - h(Y\uii_{T_{i-1}}-1) - \int_{T_{i-1}}^{S_l} \a_{Y\uii_u}\,du\,,
\en 
linking the development of the process~$h(Y\uii_t)$ with time, since, in view of~\Ref{ADB-parameter-assns-gen},
the integral is close to $\a(S_l - T_{i-1})$ if $Y\uii_{T_{i-1}}$ is large.  Since 
$\var(E'_k) =  (j + k-1)^{-2}$, it follows that, for each $l \ge 1$,  
\eq\label{ADB-M-variance}
    \var(M_l \giv \cF_{T_{i-1}}) \Eq \sum_{k=1}^l \var(E'_k) \Le 2 \Xss{i-1}^{-1}.
\en
Hence, for any $a > 0$, by Kolmogorov's inequality,
\eq\label{ADB-M-bound}
   \pr[\sup_{l \ge 0}|M_l| > a \giv \cF_{T_{i-1}}] \Le \frac{2}{a^2 \Xss{i-1}}\,.
\en
This inequality controls the fluctuations of~$(h(Y\uii_{S_l}-1) - h(Y\uii_{T_{i-1}}-1),\,l\ge1)$ around the values
defined by the integrals $(\int_{T_{i-1}}^{S_l} \a_{Y\uii_u}\,du,\,l \ge 1)$.

The next step is to bound the differences
\[
    \Delta_l \Def \int_{T_{i-1}}^{S_l} \a_{Y\uii_u}\,du - \a(S_l - T_{i-1}), \qquad l\ge 1.
\]
We have, for any $l \ge 1$,
\eqa
    |\ex\{\Delta_l \giv \cF_{T_{i-1}}\}| &=& \biggl|\sum_{j=\Xss{i-1}}^{\Xss{i-1}+l-1} j^{-1}(1 - \a/\a_j)\biggr| 
                    \Le  C_1(\Delta) \Xss{i-1}^{-\g}; \non\\
    \var(\Delta_l \giv \cF_{T_{i-1}}) &=& \sum_{j=\Xss{i-1}}^{\Xss{i-1}+l-1} j^{-2}(1 - \a/\a_j)^2 \Eq O(\Xss{i-1}^{-1-2\g}) ,
     \label{ADB-Delta-mean-var-bnd}
\ena
from which it also follows, using Kolmogorov's inequality, that
\eq\label{ADB-Delta-bound}
    \pr\bigl[\sup_{l \ge 0}|\Delta_l - \ex \Delta_l| > a \giv \cF_{T_{i-1}}\bigr] 
                   \Le \frac{C_2(\Delta)}{a^2 \Xss{i-1}^{1+2\g}}\,,
\en
for any $a > 0$, where~$C_l(\Delta)$, $l=1,2$, are constants depending on $(\a_j,\,j\ge1)$.  Hence, for any
$a > 2C_1(\Delta) X_{T_{i-1}}^{-\g}$,
\eq\label{ADB-int-bound}
    \pr\biggl[\sup_{l \ge 0} \Bigl|\int_{T_{i-1}}^{S_l} \a_{Y\uii_u}\,du - \a(S_l - T_{i-1})\Bigr| > a 
                   \giv \cF_{T_{i-1}}\biggr]
           \Le \frac{4C_2(\Delta)}{a^2 \Xss{i-1}^{1+2\g}}\,.
\en

In order to recover a bound for the probability of the event~$A_i$ from \Ref{ADB-M-bound} and~\Ref{ADB-int-bound},
we first observe that, by a simple integral comparison, for $l > k$, we have
\[
      h(l) - h(k) \Le \int_k^l x^{-1}\,dx \Eq \log l - \log k \Le h(l-1) - h(k-1),
\]
implying that $|(h(l-1)-h(k-1)) - (\log l - \log k)| \le k^{-1}$.  Then, because~$Y\uii$ is a step function,
\[
    \sup_{T_{i-1} \le t \le S_l}|h(Y\uii_{t}-1) - h(Y\uii_{T_{i-1}}-1) - \a(t - T_{i-1})|
\]
is attained at one of the points $\{S_k-,S_k,\,1\le k\le l\}$, 
and $|h(Y\uii_{S_k}) - h(Y\uii_{S_k-})| \le \Xss{i-1}^{-1}$.
Hence
\eq\label{ADB-A-prob-1}
    \sup_{t \ge t_0}|\log(Y\uii_{t}) - \log(\Xss{i-1}) - \a(t - T_{i-1})| 
     \Le \sup_{l \ge 0}|M_l + \Delta_l| + 2\Xss{i-1}^{-1}.
\en
Combining \Ref{ADB-M-bound}, \Ref{ADB-int-bound} and~\Ref{ADB-A-prob-1}, and because $X_t = Y\uii_t$ for
$T_{i-1} \le t < T_i$, it follows that, if $\Xss{i-1} \ge y_0$, where $y_0^{3\g/4} := 6\max\{C_1(\D),1\}$,
then
\[
   \pr[A_i^c \giv \cF_{T_{i-1}}] \Le \frac{2 + 4C_2(\D)}{\Xss{i-1}^{1-\g/2}} \,.
\]
This establishes the first inequality, with~$C_A:= \max\{2 + 4C_2(\D),x_0^{1-\g/2}\}$.

The inequality~(2) is much easier.  We first note that
\[
    \pr[B_i^c \giv \s(T_i \vee \cF_{T_{i}-}) \cap \{X_{T_i-}=j\}] \Eq \pr[|\log Y_j - \log p - \log j| > j^{-\cst}],
\]
where~$Y_j\sim \Pi_j$, as before.  Now, by~\Ref{ADB-parameter-assns}, 
\[
    \pr[|\log Y_j - \log p - \log j| > j^{-\cst}] \Le \pr[|\log Y_j - \log p_j - \log j| > \half j^{-\cst}]
\]
if $\half j^{-\cst} > c_2 j^{-\g}$, and, since $|\log(1+x)| \le 2|x|$ in $x \ge -1/2$,
\[
    \pr[|\log Y_j - \log p_j - \log j| > \half j^{-\cst}] \Le \pr[|(Y_j/jp_j) - 1| > \quarter j^{-\cst}]
             \Le \frac{4c_3}{p_j^2}\, j^{-\g/2},
\]
where the last inequality follows by Chebyshev's inequality, and because of~\Ref{ADB-parameter-assns}. 
Inequality~(2) now follows directly. \ep
\end{proof}

The implication of Lemma~\ref{ADB-A-and-B-probs} is that, if~$\Xss{i-1}$ is large, then $A_i \cap B_i$
occurs with high probability, in which case~$\log X_{t}$ stays close to $\a(t-T_{i-1})  + \log \Xss{i-1}$
for $T_{i-1} \le t < T_i$, and $\log \Xss{i}$ is close to $\a(T_i-T_{i-1})  + \log \Xss{i-1} - \log(1/p)$.
If this is the case for all $1\le j\le i$, then $\log X_t$ stays close to 
$\a(t-T_0) - (Z_t - Z_{T_0})\log(1/p) + \log \Xss{0}$ for $0 \le t \le T_i$.   
Now $\ex\{Z_t - Z_{T_0}\} \approx \bet(t-T_0)$,
and so $\a(t-T_0) - \ex\{Z_t - Z_{T_0}\}\log(1/p) \approx \bb(t-T_0)$, where $\bb$ 
is as in~\Ref{ADB-beta-def}.  If
$\bb > 0$, this suggests that~$\log X_t$ grows linearly with~$t$, or, equivalently, that~$X_t$
grows exponentially with~$t$, on the event 
\[
     H \Def H(j_0,t_0) := \bigcap_{i\ge1}\bigl\{ A_i(j_0,t_0) \cap B_i(j_0,t_0) \bigr\}.
\]
The remainder of the proof consists of making this heuristic precise.  

We begin by investigating the times~$T_j$ \adbb{defined in~\Ref{ADB-T-defs}}, again for the process~$(X,Z)$ starting 
a time~$t_0$ in an arbitrary 
state~$(j_0,k_0)$.  Let~$j'$ be such that, for all $j \ge j'$, 
\eq\label{ADB-jstar-def}
     \a - \bet_j\log(1/p) \ >\  \bb/2 \ >\ 0; 
\en
write  $\bet' :=  \sup_{k \ge j'}\bet_k$.
Observe that the process~$(T_j,\,j \ge 1)$ is a point process adapted to the filtration $(\cF_t,\,t \ge t_0)$, with
compensator $A(t) := \int_{t_0}^t \bet_{X_u}\,du$.  Hence, as long as $X_u \ge j'$ for all~$u$, it can be coupled 
to a point process~$(T'_j,\,j\ge1)$ 
with compensator $A'(t) := \bet'(t-t_0)$, hence a Poisson process with rate~$\bet'$, in such a way that
$T_j -T_{j-1} \ge T'_j - T'_{j-1}$ a.s.\ for all $j \ge 1$, where we take $T_0 := T_0' := t_0$.  
More precisely, taking a sequence $(U_j,\,j\ge1)$ of independent uniform $\mathrm{U}[0,1]$ random variables,
independent of everything else, define 
\eqa
    T'_j - T'_{j-1} &:=& -\frac1{\bet'}\log U_j \Eq \inf\{t \ge 0\colon \bet't \ge -\log U_j\}; \non\\
    T_j - T_{j-1}   &:=&  \inf\Bigl\{t \ge 0\colon  \int_{T_{j-1}}^t \bet_{X_u}\,du \ge -\log U_j\Bigr\},
                 \label{ADB-T-defs-2}
\ena
so that the times of the $Z$-jumps in $(X,Z)$ are defined in terms of the sequence $(U_j,\,j\ge1)$;
the times of the remaining transitions in $(X,Z)$ are defined using an independent sequence $(U'_j,\,j\ge1)$
of uniform random variables.  On the event $\Xss{j-1} \ge j'$, we then have $T_j -T_{j-1} \ge T'_j - T'_{j-1}$,
since~$X$ is non-decreasing in $[T_{j-1},T_j)$ .

For any $\h > 0$ and with~$\g$ as in~\Ref{ADB-gamma-def}, define the events
\bali
   \cE_i^{\h} &\Def \cE_i^{\h}(j_0,t_0) \Def 
           \Bigl\{\min_{1\le j\le i}\{\a(T_j'-t_0) - j(\log(1/p)+\h) + \log j_0\} > (4/\g) \log(2/\h) \Bigr\}; \non\\
    H_i &\Def H_i(j_0,t_0) \Def \bigcap_{j=1}^i \bigl(A_j(j_0,t_0) \cap B_j(j_0,t_0)\bigr). \label{ADB-Ei-def}
\end{align}
The next step is to show that, on the event $\cE_i^{\h}(j_0,t_0) \cap H_i(j_0,t_0)$, the values of $\Xss{i} = X_{T_i}$
can be bounded below in terms of the times~$T'_i$.

\begin{lemma}\label{ADB-first-growth}
If $\bb > 0$, choose any $\h > 0$ such that
$(2/\h)^{4/\g} > j'$, where~$j'$ is as defined in~\Ref{ADB-jstar-def}.
Then, for any $j_0 > (2/\h)^{4/\g}$, on the event $\cE_i^{\h}(j_0,t_0) \cap H_i(j_0,t_0)$,  it follows that
\eq\label{ADB-growth-1}
    2\Xss{j}^{-\cst}\ <\ \h  \andy \log \Xss{j}\ \ge\ \a(T_j' - t_0) - j(\log(1/p)+\h) + \log j_0
\en
for all $0\le j\le i$.
\end{lemma}

\begin{proof}
 The proof runs by induction on~$j$.  By assumption, $2\Xss{0}^{-\cst} = 2j_0^{-\cst} < \h$, and
the second inequality is trivial for~$j=0$.  If~\Ref{ADB-growth-1}, with $j-1$ for~$j$, is true for some~$j$, then, 
from the definitions~\Ref{ADB-AB-defs} of $A_j(j_0,t_0)$ and~$B_j(j_0,t_0)$, 
\bals
      \log \Xss{j} - \log\Xss{j-1}\ &\ge\ \a(T_j - T_{j-1}) - \log(1/p) - 2\Xss{j-1}^{-\cst} \\
                   &\ge\ \a(T_j' - T_{j-1}') - (\log(1/p) + \h);
\end{align*}
here, we have used $2\Xss{j-1}^{-\cst} < \h$, so that $\Xss{j-1} \ge j'$; and then also the fact that~$X$ is non-decreasing 
on $[T_{j-1},T_j)$, so that 
$\bet_{X_u} \le \bet'$ for all $T_{j-1} \le u < T_j$.  Adding this inequality to the second induction hypothesis gives 
\[
    \log \Xss{j}\ \ge\ \a(T_{j}' - t_0) - j(\log(1/p)+\h) + \log j_0 \ \ > \ (4/\g)\log (2/\h),
\]
in view of~$\cE_i^{\h}(j_0,t_0)$, implying that $2\Xss{j}^{-\cst} < \h$ also. \ep
\end{proof}

The lower bound established in Lemma~\ref{ADB-first-growth} is useful for establishing that the values~$\Xss{j}$ grow
geometrically fast with~$j$ for $1\le j\le i$, on the event $\cE_i^{\h}(j_0,t_0) \cap H_i(j_0,t_0)$, provided that 
the distribution of the times $(T_j',\, j\ge 1)$ can be controlled.  This is made possible by the next lemma.

\begin{lemma}\label{ADB-ex-min}
 Let $(E_i,\,i\ge1)$ be independent standard exponential random variables, and let $V_i := \sum_{j=1}^i E_j$.

\smallskip
{\rm (i)}. For any $0 < \f < 1$, let~$u(\f)$ be the positive solution to $u\f = \log(1+u)$.  Then, for any $a > 0$,
\[
             \pr[\min_{j \ge 1}\{V_j - j\f\} \le - a] \Le e^{-au(\f)}.
\]

{\rm (ii)}. For any $0 < \f < 1$ and $c,a > 0$, there exist $k = k(\f,c) > 0$ and $0 < \ps = \ps(\f,c) < 1$ such that
\[
    \sum_{i\ge i_0} \ex\bigl\{ \min\bigl(1,\,e^{-c(V_i - i\f + a)}\bigr) \bigr\} 
            \Le k \ps^{\max\{i_0,a\}}  ,      \quad i_0 \ge 1.
\]
Hence, in particular,
\[
    \sum_{i\ge 1}  \min\bigl(1,\,e^{-c(V_i - i\f + a)}\bigr) \ <\ \infty\ \mbox{a.s.}
\]
\end{lemma}

\begin{proof}
For part~(i), it is immediate that the sequence $M_j^u := e^{-uV_j + j\log(1+u)}$ is a martingale for
any $u > -1$.  Writing $\st(a) := \inf\{j \ge 0\colon V_j - j\f \le -a\}$, it thus follows that,
for any $n \ge 1$,
\eqs
  1 &=& \ex\{M^{u(\f)}_{\st(a)\wedge n}\} \Eq \ex\bigl\{ e^{-u(\f)(V_{\st(a)\wedge n} - (\st(a)\wedge n)\f)} \bigr\}
              \ \ge\ \pr[\st(a) \le n]e^{au(\f)},
\ens
from which part~(i) follows.

\adbb{
For part~(ii), observe first that, for any $\f' \in (\f,1)$, we have
\eqa
   \ex\bigl\{ \min\bigl(1,\,e^{-c(V_i - i\f + a)}\bigr) \bigr\} &\le&
       \pr[V_i \le i\f' - a/2] + \ex\bigl\{e^{-c(V_i - i\f + a)}I[V_i > i\f' - a/2]\bigr\} \non\\
        &\le& \pr[V_i \le i\f' - a/2] + e^{-ci(\f' - \f) - ca/2}.  \label{ADB-exm-1}
\ena
If $i\f' \le a/2$, $\pr[V_i \le i\f' - a/2] = 0$.  For $i\f' > a/2$, take $u_2(\f') := \half u(\f')$, for which
\[
     v(\f') \Def \log(1 + u_2(\f')) - \f' u_2(\f') \ >\ 0,
\]
beacuse the function $\log(1+u)$ is strictly concave, and $\log(1+u) = u\f'$ for $u = 0$ and~$u(\f')$.  Then
\eqs
    1 &=& \ex\{M_i^{u_2(\f')}\} \ \ge\ \ex\bigl\{\exp\{-u_2(\f')V_i + i\log(1 + u_2(\f'))\}I[V_i \le i\f' - a/2]\bigr\} \\
      &\ge& e^{au_2(\f')/2 + i v(\f')}\pr[V_i \le i\f' - a/2].
\ens
Hence
\eq\label{ADB-exm-3}
    \sum_{i \ge i_0}\pr[V_i \le i\f' - a/2] \Le e^{-a u_2(\f')/2}\,\frac{e^{-i_0 v(\f')}}{1 - e^{-v(\f')}}\,.
\en
Since also
\[
    \sum_{i \ge i_0} e^{-ci(\f' - \f) - ca/2} \Le e^{-ca/2}\,\frac{e^{-i_0 c(\f'-\f)}}{1 - e^{-c(\f'-\f)}}\,,
\]
the second part follows from~\Ref{ADB-exm-1}.  \ep
}
\end{proof}

The next step is to use Lemmas \ref{ADB-A-and-B-probs}--\ref{ADB-ex-min} to conclude that,
by choosing~$j_0$ large enough, the probability $\pr[H(j_0,t_0) \cap \cE^{\h}(j_0,t_0)]$ can be made
arbitrarily close to~$1$;  here, $\cE^{\h}(j_0,t_0) := \bigcap_{i\ge1} \cE_i^{\h}(j_0,t_0)$.

\begin{lemma}\label{ADB-prob-HE}
 Assume that $\bb := \a - \bet\log(1/p) > 0$, and choose any $0 < \th < 1$ and $\h > 0$ 
such that
\eq\label{ADB-HE-1}
    \h \Le \frac{\th\bb'}{\bet'} \andy \Bigl(\frac2{\h}\Bigr)^{4/\g} \ >\ j',
\en
where $\bb' := \a - \bet'\log(1/p) > 0$, and $\bet'$ and~$j'$ are as in~\Ref{ADB-jstar-def}. Then
\[
     \pr[H(j_0,t_0) \cap \cE^{\h}(j_0,t_0)] \ \ge\ 1 - C_0 j_0^{-\d_0},
\]
for $0 < \d_0,C_0 < \infty$ depending only on the sequences $(\a_j,\bet_j,p_j,\s^2_j)$  and~$\th$.
\end{lemma}

\begin{proof}
 First, we show that $\pr[(\cE^{\za}(j_0,t_0))^c]$ is small if~$j_0$ is large, 
for all~$\za \le \h$. 
Writing $\f_{\za} := \bet'(\log(1/p) + \za)/\a \le \f_{\h} < 1$,
the process $(\a(T_j'-t_0) - j(\log(1/p)+\za),\,j\ge0)$ is that of the partial sums of independent random variables
$(\a/\bet')(E_j - \f_{\za})$, where $(E_j,\,j\ge1)$ is a sequence of independent exponential random variables
with mean~$1$. Hence, by Lemma~\ref{ADB-ex-min}~(i), for any $a > 0$,
\[
     \pr[M_\za \le -a] \Le e^{-au(\f_{\za})\bet'/\a}, \quad\mbox{where}\quad 
                M_\za \Def \min_{j\ge1}\{\a(T_j'-t_0) - j(\log(1/p)+ \za)\} .
\]
Hence
\bali
    \pr[(\cE^{\za}(j_0,t_0))^c] \Eq \pr[M_\za \le (4/\g) \log(2/\za) - \log j_0] 
             &\Le  \biggl\{\Bigl(\frac{\za}2 \Bigr)^{4/\g} j_0\biggr\}^{-u(\f_{\za})\bet'/\a} . \label{ADB-Ec-bound}
\end{align}

Next, we express $\pr[H(j_0,t_0)^c \cap \cE^{\h}(j_0,t_0)]$ in the form
\bali
      \pr[H^c \cap \cE^\h] &\Eq \pr\Bigl[ \Bigl(\bigcap_{i\ge1} H_i\Bigr)^c \cap \cE^\h\Bigr] \non\\
        &\Eq \pr[(A_1 \cap B_1)^c \cap \cE^\h] + \sum_{i\ge1}\pr[H_i \cap (A_{i+1} \cap B_{i+1})^c \cap \cE^\h] .
              \label{ADB-Hc}
\end{align}
Now, because $\cE^\h \subset \cE^\h_i$ for each~$i$, we have
\bali
        \pr[H_i \cap (A_{i+1} \cap B_{i+1})^c \cap \cE^\h] &\Le \pr[H_i \cap (A_{i+1} \cap B_{i+1})^c \cap \cE^\h_i] \non\\
               & \Eq \ex\{\pr[H_i \cap (A_{i+1} \cap B_{i+1})^c \cap \cE^\h_i \giv \cF_{T_i}]\}, \label{ADB-Hc-1}
\end{align}
and it is immediate that 
\bali
    \pr[H_i \cap (A_{i+1} \cap B_{i+1})^c \cap \cE^\h_i \giv \cF_{T_i}] &\Eq
              \pr[(A_{i+1} \cap B_{i+1})^c \giv \cF_{T_i}]I[H_i \cap \cE^\h_i] . \label{ADB-Hc-2}
\end{align}
From Lemma~\ref{ADB-A-and-B-probs} and both inequalities in Lemma~\ref{ADB-first-growth}, it thus follows that
\bali
    \pr[H_i &\cap (A_{i+1} \cap B_{i+1})^c \cap \cE^\h_i \giv \cF_{T_i}] 
              \Le (C_A + C_B)\Xss{i}^{-(\g/2)} I[H_i \cap \cE^\h_i] \non\\
           &\Le (C_A + C_B)
                \min\bigl\{ 1, \exp\{-(\g/2)(\a(T_i' - t_0) - i(\log(1/p)+\h) + \log j_0)\} \bigr\}, 
                       \label{ADB-Hc-3}
\end{align}
where we have used $j' \ge 1$.  Combining \Ref{ADB-Hc-1} and~\Ref{ADB-Hc-3}, it follows that
\bali
    \pr[H_i \cap (A_{i+1} \cap B_{i+1})^c \cap \cE^\h] 
              &\Le (C_A + C_B)\ex\bigl\{1 \wedge e^{-c(V_i - i\f_\h + a_0)}\Bigr\}, 
               \label{ADB-Hc-4}
\end{align}
where $V_i$ and~$\f_a$ are as before, 
\[
    c \Def \frac{\a\g}{2\bet'} \ >\ 0 \andy  a_0 \Def \frac{\bet'\log j_0}\a\,.
\]
Applying Lemma~\ref{ADB-ex-min}~(ii), it follows from \Ref{ADB-Hc} and~\Ref{ADB-Hc-4} that
\eq\label{ADB-Hc-5}
     \pr[H^c \cap \cE^\h] \Le (C_A + C_B) k(\f_\h,c) \{\ps(\f_\h,c)\}^{a_0};
\en
taking $\za =\h$ in~\Ref{ADB-Ec-bound}, together with~\Ref{ADB-Hc-5} and the definition of~$a_0$, proves the lemma. \ep
\end{proof}

With these preparations, it is now possible to prove Theorem~\ref{aslimit}.

\medskip
\begin{prooff} {\bf of Theorem~\ref{aslimit}}.
 For any $j \in \nat$, let $T_0 = \t(j) := \inf\{t > 0\colon X_t = j\} \le \infty$.  Then, defining
$\cX_j := \{\t(j) < \infty\}$, we have $\pr[H(j,\t(j)) \cap \cE^{\h}(j,\t(j)) \giv \cX_j] \ \ge\ 1 - C_0 j^{-\d_0}$, by
Lemma~\ref{ADB-prob-HE} and the strong Markov property.  On the event $H(j,\t(j)) \cap \cE^{\h}(j,\t(j))$,
for any $r \ge 0$, we have
\eq\label{ADB-pf-of-th-1}
    \sup_{t \ge T_r}|\log X_t - \a(t-T_r) + (Z_t - Z_{T_r})\log(1/p) - \log\Xss{r}| \Le 2\sum_{l\ge r}\Xss{l}^{-\cst},
\en
and hence, \adbb{recalling that $W_t := e^{-\alpha t}p^{-Z_t} X_t$, it follows that}
\eq\label{ADB-pf-of-th-1.2}
    \adbb{\sup_{s,t \ge T_r}|\log W_t - \log W_s| \Le 4\sum_{l\ge r}\Xss{l}^{-\cst}.}
\en
Then, by Lemma~\ref{ADB-first-growth}, on the event $H(j,\t(j)) \cap \cE^{\h}(j,\t(j))$, we have
\eq\label{ADB-pf-of-th-1.5}
   \Xss{i}\ \ge\ \max\bigl\{1,\,j e^{\a(T_i' - T_0) - i(\log(1/p)+\h)} \bigr\} \quad\mbox{for all}\ i \ge 0,
\en
and this, in turn, from Lemma~\ref{ADB-ex-min}, implies that $\sum_{i\ge 1}\Xss{i}^{-\cst} < \infty$ a.s., \adbb{and hence 
that, on the event $H(j,\t(j)) \cap \cE^{\h}(j,\t(j))$,}
\eq\label{ADB-pf-of-th-2}
   \adbb{\lim_{r\to\infty} \sum_{i\ge r}\Xss{i}^{-\cst} \Eq 0\ \ \mbox{a.s.} } 
\en
Thus, on the event 
\[
   A(j) \Def \cX_j \cap H(j,\t(j)) \cap \cE^{\h}(j,\t(j)),
\]
we conclude \adbb{from \Ref{ADB-pf-of-th-1.2} and~\Ref{ADB-pf-of-th-2}, using the Cauchy criterion,} that
$\log W_t$ converges a.s.\ to a finite limit, implying that~$W_t$ converges a.s.\ to a strictly positive limit~$W$.

Now, for $X_0 = x_0$, $\cX = \bigcap_{j > x_0} \cX_j$. Writing $\pr_j[\cdot]$
for $\pr[\cdot\giv X_0 = j]$, this implies that
$\pr_{x_0}[\cX] = \lim_{j\to\infty}\pr_{x_0}[\cX_j]$, and, as above, 
$A(j) \subset \cX_j$ and, using Lemma~\ref{ADB-prob-HE},
\[
  0 \Le \pr_{x_0}[\cX_j] - \pr_{x_0}[A(j)] \Le \pr_{x_0}[\cX_j] C_0 j^{-\d_0} \ \to\ 0
            \quad\mbox{as}\ j \to \infty .
\]
Hence, if $C := \{W_t \mbox{ converges to a positive limit}\}$, we have $A(j) \subset C $ for each~$j$,
and 
\[
    \limsup_{j\to\infty}\pr_{x_0}[A(j)]  \Eq \lim_{j\to\infty}\pr_{x_0}[\cX_j] \Eq \pr_{x_0}[\cX].
\]
Hence $\pr_{x_0}[C] \ge \pr_{x_0}[\cX]$.  However, $\cX^c = \{X_t = 0 \mbox{ eventually}\}$, on which event
$W_t \to 0$ a.s.  Hence $C = \cX$ a.s., and the
theorem is proved.  \ep
\end{prooff}

\section{The path to infinity: the central limit theorem}\label{ADB-sect-CLT}
\setcounter{equation}{0}

\adbb{
In view of Theorem~\ref{aslimit}, and its interpretation in~\Ref{ADB-heuristic-2}, the main fluctuations
in $\log X_t$, for large~$t$, are governed by those of a Poisson process.  In keeping with this interpretation,
we now deduce a central limit theorem for~$\log X_t$, which is based on the central limit theorem for the Poisson
distribution.  As a consequence, for large~$t$, that part of the distribution of~$X_t$ that is not the mass on zero
has an approximately log-normal distribution.
}

\begin{theorem}\label{CLT}
 Under the assumptions of Theorem~\ref{aslimit}, for any $x_0 \in \nat$, the process~$X$ with $X_0 = x_0$
satisfies
\[
     \lti \pr_{x_0}[t^{-1/2}\{\log X_t - \bb t\} \ge y] \Eq  \pr_{x_0}[\cX](1 - \Phi(y/v)),
\]
where $v^2 := \bet\{\log(1/p)\}^2$ and~$\bb$ is as in~\Ref{ADB-beta-def}.
\end{theorem}

\begin{proof}
From the definition of~$W_t$, we have
\[
    \log X_t \Eq \a t - Z_t \log(1/p) - \log W_t,
\]
so that, for any $0 < s < t$,
\begin{align}
   \log X_t - \bb t &\Eq - (Z_t - \bet t)\log(1/p) - \log W_t \label{ADB-CLT-1}\\
     &\Eq - (Z_t - Z_s - \bet(t-s))\log(1/p) - (Z_s - \bet s)\log(1/p) - \log W_t.
         \non
\end{align}
The remainder of the proof consists of showing that, for suitable choice of~$s$, the final two terms are 
asymptotically negligible, and that $Z_t - Z_s$ is approximately Poisson distributed with mean $\bet(t-s)$.

\medskip
\nin{\bf The last two terms in~\Ref{ADB-CLT-1}}\\
We start by showing that, if $t \gg s \gg 1$, the final two terms of~\Ref{ADB-CLT-1} are small \adbb{compared to~$\sqrt t$} 
on a `good event' that has probability close enough to $\pr_{x_0}[\cX]$.  In preparation,
choose any $0 < \e < \half\bet$.  With~$j'$ is as defined in~\Ref{ADB-jstar-def},   
choose~$j_0 := j_0(\e) \ge j'$ such that $c_1 j_0^{-\g_1} < \e$. 
Then choose~$j_1 := j_1(\e) \ge (j_0(\e))^2$ so that, with $\cX_j := \{\t(j) < \infty\}$,  
\eq\label{ADB-X-cvgce-prob}
   0 \ <\ \pr_{x_0}[\cX] - \pr_{x_0}[\cX_{j_1}] \ <\ \e,
\en 
this is possible, since $\cX = \bigcap_{j > x_0}\cX_j$.
Choose 
\bali
   s &\Def s(\e)\ >\ 0 \quad \mbox{such that} \quad \pr_{x_0}[\t(j_1(\e)) > s \giv \cX_{j_1}]\ <\ \e; \non\\
   z &\Def z(\e)\ >\ 0 \quad \mbox{such that} \quad \pr_{x_0}[Z_{s(\e)} > z \giv \t(j_1) \le s(\e)]\ <\ \e;
                           \label{ADB-s,z,w-defs} \\
   w &\Def w(\e)\ >\ 0 \quad \mbox{such that} \quad \pr_{x_0}[|\log W_{\t(j_1)}| > w \giv \t(j_1) \le s(\e)]\ <\ \e, \non 
\end{align}
\adbb{noting that the quantities $j_0,j_1,s,z$ and~$w$ depend on the choice of~$\e$, but are then fixed.
Let $A_\e\ui := \{\t(j_1(\e)) \le s(\e)\}$, and
define part of the `good event' by
\eq\label{ADB-E1-def}
   E\ui_\e \Def A_\e\ui \cap \{Z_{s(\e))} \le z(\e)\} \cap \{|\log W_{\t(j_1(\e))}| \le w(\e)\};
\en
on~$E\ui_\e$,} the values of the process up to the time of first hitting the state~$j_1(\e)$ are not extreme.

Take $T_0 := T_0(\e) := \t(j_1(\e))$, and let the subsequent jumps of~$Z$ be denoted by $(T_i,\,i\ge1)$; as before, 
note that~$T_0$ itself need not be a jump time of~$Z$.  Couple them to
the jump times $(T_i',\,i\ge1)$ of a Poisson process~$Z'$ of rate $\bet'$ as in~\Ref{ADB-T-defs-2}, and,
fixing  $\h > 0$ to satisfy~\Ref{ADB-HE-1}, define
\eqa
    \hcE_i^{\h} &:=& \hcE_i^{\h}(j_1,T_0) \Def 
           \Bigl\{\min_{1\le j\le i}\{\a(T_j'-T_0) - j(\log(1/p)+3\h/2)\} >  - \half\log j_1 \Bigr\}; \non\\
    \hcE^{\h} &:=& \hcE^{\h}(j_1,T_0) \Def \bigcap_{i \ge 1} \hcE_i^{\h}(j_1,T_0). \label{ADB-E-hat-def}
\ena
The event 
\eq\label{ADB-E2-def}
   \adbb{E_\e\ut} \Def  H(j_1(\e),\t(j_1(\e))) \cap \cE^{\h}(j_1(\e),\t(j_1(\e))) \cap \hcE^{\h}(j_1(\e),\t(j_1(\e))), 
\en
is to be the second element of the `good event'.
On this event, using~\Ref{ADB-pf-of-th-1.5} and the definition of~$\hcE^{\h}$, it follows that
\eq\label{ADB-pf-of-CLT-1.5}
   \Xss{i}\ \ge\ j_0 e^{i\h/2}  \quad\mbox{for all}\ i \ge 0.
\en
Thus, on the event $\cX_{j_1(\e)} \cap E\ut_\e$, from~\Ref{ADB-pf-of-th-1} and~\Ref{ADB-pf-of-CLT-1.5},
\bali
   \sup_{t \ge \t(j_1(\e))}|\log X_t - &\log j_1(\e) - \a(t-\t(j_1(\e))) + (Z_t - Z_{\t(j_1(\e))})\log(1/p) | \non\\
        &\Eq   \sup_{t \ge \t(j_1(\e))}|\log W_t - \log W_{\t(j_1(\e))}|  
                  \Le   \frac{2}{(j_0(\e))^{\cst}(1-e^{-\h\g/8})}\,,  \label{ADB-W-bound}
\end{align}
\adbb{
implying that, on the event $\cX_{j_1(\e)} \cap E\ui_\e \cap E\ut_\e$,
\eq\label{ADB-logW-bnd}
       \sup_{t \ge \t(j_1(\e))}|\log W_t| \Le w(\e) + \frac{2}{(j_0(\e))^{\cst}(1-e^{-\h\g/8})}\,.
\en
Hence, on the same event, for all $t \ge s(\e)$,
\eq\label{ADB-extras-bnd}
       |Z_{s(\e)} - \bet s(\e)|\log(1/p) + |\log W_t| \Le 
               (z(\e)+ \bet s(\e))\log(1/p) + w(\e) + \frac{2}{(j_0(\e))^{\cst}(1-e^{-\h\g/8})}\,;
\en
thus, as $t\to\infty$ and with $s = s(\e)$, the contribution to~\Ref{ADB-CLT-1} from the last two terms is asymptotically
negligible compared to~$\sqrt t$.}

\adbb{
Note also that, from~\Ref{ADB-W-bound}, we have
\bali
   \sup_{s,t \ge \t(j_1(\e))}|\log X_t - &\log X_s - \a(t-s) + (Z_t - Z_s)\log(1/p) | \non\\
        &\Le   2\sup_{t \ge \t(j_1(\e))}|\log W_t - \log W_{\t(j_1(\e))}|  
                  \Le   \frac{4}{(j_0(\e))^{\cst}(1-e^{-\h\g/8})}\,.  \label{ADB-W-bound-2}
\end{align}
}
This shows that the development of $\log X_t - \a t$ after the time~$\t(j_1(\e))$ is mirrored in the development
of~$-Z_t$. 
Furthermore,  on the event~$\cX_{j_1(\e)} \cap E_\e\ut$ and from~\Ref{ADB-parameter-assns},
for $i \ge 1$ and for $T_{i-1} \le t < T_i$,
\eq\label{ADB-q-bound-2}
     |\bet_{X_t} - \bet| \Le c_1 (j_0(\e) e^{(i-1)\h/2})^{-\g_1} \Le \e e^{-(i-1)\h\g_1/2}\ =:\ \e_\bet(i).
\en

\medskip
\nin{\bf For large times, $Z$ is close to a Poisson process with constant rate}\\
The next step is to show that, to approximate the development of $\log X_t - \a t$ after time~$\t(j_1(\e))$,
the process~$Z$ can be replaced by a Poisson process of rate $\bet$, with only small error.
\adbb{Suppressing the dependence on~$\e$ in the notation, where possible,
we set $T_0 := \t(j_1)$, as before, and set $\tT_0 := T_0$.}
We then use the independent uniform random variables $(U_i,\,i\ge1)$ introduced in~\Ref{ADB-T-defs-2} to define the jump times
$(\tT_i,\,i\ge1)$ of a further Poisson process~$\tZ$, of rate $\bet$, on the time interval $(T_0,\infty)$, coupled to 
the jump times $(T_i,\,i\ge1)$ of~$Z$ on this time interval.  The coupling is achieved by defining
$\tT_i - \tT_{i-1} := -\frac1{\bet}\log U_i$, with~$T_i - T_{i-1}$ defined as in~\Ref{ADB-T-defs-2}.  It then follows
from~\Ref{ADB-q-bound-2} that, on the event~$\cX_{j_1} \cap E\ut(j_1)$,
\[
     |(\tT_i - \tT_{i-1}) - (T_i - T_{i-1})| \Le \frac{\e_\bet(i)}{\bet-\e_\bet(i)}\,(\tT_i - \tT_{i-1})
            \Le \frac{2\e_\bet(i)}{\bet}\,(\tT_i - \tT_{i-1}),
\]
because $\e < \half\bet$, and hence that
\[
     |\tT_i - T_i| \Le \sum_{l=1}^i \frac{2\e_\bet(l)}{\bet^2}\,(-\log U_l).
\]
Because $\sum_{i\ge1} \e_\bet(i) < \infty$, it follows that
\[
      T^* \Def \sum_{l\ge 1} \frac{2\e_\bet(l)}{\bet^2}\,(-\log U_l)
\]
is an a.s.\ finite random variable, and that
\[
              \sup_{i \ge 1}|\tT_i - T_i| \Le T^*\quad \mathrm{a.s.};
\]
indeed, 
\eq\label{ADB-T-star-mean}
    \ex T^* \Eq \sum_{l\ge 1} \frac{2\e_\bet(l)}{\bet^2} \Le \frac{2\e}{\bet^2(1 - e^{-\h\g/2})}\ =:\ \e t^*,
\en 
also.
Thus, on the event~$\cX_{j_1} \cap E_\e\ut$, for all $t \ge T^* + T_0 = T^* + \t(j_1)$, we have
\eq\label{ADB-Z-Poisson-approx}
    \tZ_{t - T^*} - \tZ_{\t(j_1)} \Le Z_t - Z_{\t(j_1)} \Le  \tZ_{t + T^*} - \tZ_{\t(j_1)}.
\en
This in turn implies that
\eq\label{ADB-Z-Poisson-approx-2}
    |(Z_t - Z_{\t(j_1)}) - (\tZ_{t} - \tZ_{\t(j_1)})| \Le \tZ_{t + T^*} - \tZ_{t - T^*} \ =:\ \tZ^*_t,
\en
bounding the error involved at time~$t$, if $Z$ is replaced by~$\tZ$.
Note that $\pr[T^* > 2t^*] \le \half \e$, by Markov's inequality, and that, on $\{T^* \le 2t^*\}$,
\[
  \tZ^*_t \Le \tZ_{t+2t^*} - \tZ_{t - 2t^*}\ \sim\ \Po(4t^*\bet) \quad\mbox{if}\quad t \ge 2t^* + \t(j_1).  
\]
Hence, choosing $\tz := \tz(\e)$ in such a
way that $\Po(4t^*\bet)\{[\tz,\infty)\} \le \half\e$, it follows that, for all $t \ge  s(\e) + 2t^*$,
\eq\label{ADB-tZs-prob}
    \pr_{x_0}[\tZ^*_t > \tz \giv A\ui_\e ]\ <\ \e,
\en
quantifying the error~\Ref{ADB-Z-Poisson-approx-2}.

\medskip
\nin{\bf Approximation using Poisson probabilities}\\
Let $\hZ_t := -(\tZ_t - \bet t)$ denote minus the centred version of~$\tZ_t$.
Then, for $t > s(\e)$ and on the `good event' 
\eq\label{ADB-At-def}
            A_{t,\e} \Def E\ui_\e \cap E\ut_\e \cap \{\tZ^*_t\le \tz(\e)\}, 
\en
we can use \Ref{ADB-s,z,w-defs},  \Ref{ADB-Z-Poisson-approx-2} and~\Ref{ADB-W-bound-2} 
to give the following summary of the approximation of $\log X_t - \bb t$  
\adbb{using the centred Poisson random variable} $\hZ_t - \hZ_{s(\e)}$:
\eqa
    \lefteqn{|(\log X_t - \bb t) - (\hZ_t - \hZ_{s(\e)})\log(1/p)|\,I[A_{t,\e}]} \label{ADB-main-split}\\
          &&\Le  \{z(\e) + s(\e)\bet + \tz(\e)\}\log(1/p) + w(\e) + \frac{4}{j_0(\e)^{\cst}(1-e^{-\h/2})}
            \ =:\ \ue(\e).  \non
\ena
Note that the `natural' approximation to $\log X_t - \bb t$ would be $\hZ_t \log(1/p)$, but for~\Ref{ADB-V-dist} 
below we need the fact that, \adbb{conditional on any event in~$\cF_{s(\e)} \cap A\ui_\e$, 
the increment $\hZ_t - \hZ_{s(\e)}$ has the centred version of the Poisson distribution with mean $\bet(t - s(\e))$;
using this increment in the approximation} leads to the term~$\{z(\e) + s(\e)\bet\}\log(1/p)$ in
the error bound~$\ue(\e)$.  

Using~\Ref{ADB-main-split}, it follows that, for any $y \in \re$, we have
\bals
  \bigl( \{(\hZ_t - \hZ_{s(\e)})&\log(1/p) \ge y\sqrt t + \ue(\e)\} \cap A\ui_\e \bigr) \setminus (A_{t,\e}^c \cap A\ui_\e) \\
      &\subset\ \{\log X_t - \bb t \ge y\sqrt t\}\cap A_{t,\e} \\ 
      &\subset\ \{(\hZ_t - \hZ_{s(\e)})\log(1/p) \ge y\sqrt t - \ue(\e)\} \cap A\ui_\e ,
\end{align*}
and hence that
\bali
     \pr_{x_0}[\{(\hZ_t - \hZ_{s(\e)})&\log(1/p) \ge y\sqrt t + \ue(\e)\} \cap A\ui_\e] 
                    - \pr_{x_0}[A_{t,\e}^c \giv A\ui_\e] \non\\
     &\Le \pr_{x_0}[\{\log X_t - \bb t \ge y\sqrt t\} \cap A\ui_\e] \label{ADB-log-to-V}\\
     &\Le \pr_{x_0}[\{(\hZ_t - \hZ_{s(\e)})\log(1/p) \ge y\sqrt t - \ue(\e)\} \cap A\ui_\e] . \non
\end{align}
\adbb{
Because the increment $\hZ_t - \hZ_{s(\e)}$ has the Poisson distribution $\Po\bigl(\bet(t-s(\e))\bigr)$,
conditional on any event in $\cF_{s(\e)} \cap A\ui_{\e}$, it also follows that, for any $x \in \re$,
\eq\label{ADB-V-dist}
    \pr_{x_0}[\{(\hZ_t - \hZ_{s(\e)})\log(1/p) \ge x\} \cap A\ui_\e] 
                 \Eq \pr_{x_0}[A\ui_\e]\Po\bigl(\bet(t-s(\e))\bigr)[\bet(t-s(\e)) + x/\log(1/p),\infty),
\en
and this can be used on both the left and the right hand sides of~\Ref{ADB-log-to-V}, with only slightly
different values of~$x$.
}

\medskip
\nin{\bf A first normal approximation}\\
We now convert the inequalities~\Ref{ADB-log-to-V} into an explicit normal approximation to the probability
$\pr_{x_0}[\{\log X_t - \bb t \ge y\sqrt t\} \cap \cX]$.
First, by the Berry--Esseen theorem, if $t > 2s(\e)$, we have
\[
    |\pr[(\hZ_t - \hZ_{s(\e)})\log(1/p) \ge x\sqrt{t - s(\e)}] - (1 - \Phi(x/v))| 
                   \Le \frac{2C_{\bet}}{\sqrt t}\,,\quad x \in \re,
\]
where~$C_{\bet}$ is a finite constant and $v^2$ is as in the statement of the theorem.  Thus
\[
     \pr[(\hZ_t - \hZ_{s(\e)})\log(1/p) \ge y\sqrt t - \ue(\e)] \Le 
              1-\Phi\Bigl(\frac{y\sqrt t - \ue(\e)}{v\sqrt{t - s(\e)}}\Bigr) + \frac{2C_{\bet}}{\sqrt t}\,.
\]
Then, from the properties of the normal distribution, and using $t > 2s(\e)$,
\[
     0 \ <\ \Phi\Bigl(\frac{y\sqrt t}{v\sqrt{t - s(\e)}}\Bigr) - \Phi\Bigl(\frac{y\sqrt t - \ue(\e)}{v\sqrt{t - s(\e)}}\Bigr) 
              \Le \frac{\ue(\e)}{v\sqrt{\pi t}}\,, 
\]
and 
\[
       \Bigl|\Phi\Bigl(\frac{y\sqrt t}{v\sqrt{t - s(\e)}}\Bigr) - \Phi\Bigl(\frac{y}{v}\Bigr)\Bigr|
              \Le \frac{|y|e^{-(y/v)^2/2}}{v\sqrt{2\pi}}\,\Bigl\{\sqrt{\frac{t}{t-s(\e)}} - 1\Bigr\}
              \Le \frac{1}{\sqrt{2\pi e}}\, \frac{s(\e)}t \,,
\]
where we have used $0 \le (1-x)^{-1/2} - 1 \le x$ in $0 \le x \le 1/2$.
Collecting these bounds, we have shown that
\eq\label{ADB-CLT-5}
    \pr[(\hZ_t - \hZ_{s(\e)})\log(1/p) \ge y\sqrt t - \ue(\e)] \Le (1-\Phi(y/v)) + \d(t,\e),
\en
where
\eq\label{ADB-CLT-6}
    \d(t,\e) \Def  \frac{2C_{\bet}}{\sqrt t} + \frac{\ue(\e)}{v\sqrt{\pi t}} + \frac{1}{\sqrt{2\pi e}}\, \frac{s(\e)}t\,;
\en
note that $\lti \d(t,\e) = 0$ for each $\e > 0$. An entirely similar argument shows that
\eq\label{ADB-CLT-7}
    \pr[(\hZ_t - \hZ_{s(\e)})\log(1/p) \ge y\sqrt t + \ue(\e)] \ \ge\ (1-\Phi(y/v)) - \d(t,\e).
\en
Hence, in view of~\Ref{ADB-log-to-V}, we have
\eq\label{ADB-CLT-8}
   |\pr_{x_0}[\{\log X_t - \bb t \ge y\sqrt t\} \cap A_{t,\e}] - (1-\Phi(y/v))\pr_{x_0}[A\ui_\e]| 
            \Le \pr_{x_0}[A_{t,\e}^c \giv A\ui_\e] + \d(t,\e).
\en

\medskip
\nin{\bf The final normal approximation}\\
It now remains to tidy up, replacing $A_{t,\e}$ and~$A\ui_\e$ by~$\cX$ in~\Ref{ADB-CLT-8}, and evaluating
the various error terms.  First, observe that
\[
     \cX \ \subset\ \cX_{j_1(\e)} \andy A\ui_\e \Eq \{\t(j_1(\e)) \le s(\e)\} \ \subset\ \cX_{j_1(\e)},
\]
so that, for any event~$B$,
\eqa
      \lefteqn{|\pr_{x_0}[B \cap \cX] - \pr_{x_0}[B \cap A\ui_\e]|} \non\\ 
        &&\Le  (\pr_{x_0}[\cX_{j_1(\e)}] - \pr_{x_0}[\cX]) 
                    + (\pr_{x_0}[\cX_{j_1(\e)}] - \pr_{x_0}[\t(j_1(\e)) < s(\e)]) 
         \ <\   2\e,  \label{ADB-A1-B-prob}
\ena
where the last inequality follows from~\Ref{ADB-X-cvgce-prob} and from the definition of~$s(\e)$.  
In particular, taking~$B$ to be the certain event, this implies that
\eq\label{ADB-A1-prob}
      |\pr_{x_0}[\cX] - \pr_{x_0}[A\ui_\e]| \ <\ 2\e.
\en
Thus, again using~\Ref{ADB-A1-B-prob},
\bals
     |\pr_{x_0}[\{&\log X_t - \bb t \ge y\sqrt t\} \cap A_{t,\e}] 
                   - \pr_{x_0}[\{\log X_t - \bb t \ge y\sqrt t\} \cap \cX]| \\
           &\Le |\pr_{x_0}[\{\log X_t - \bb t \ge y\sqrt t\} \cap A\ui_{\e}] - 
               \pr_{x_0}[\{\log X_t - \bb t \ge y\sqrt t\} \cap \cX]| \\
           &\hskip2cm                             + \pr_{x_0}[A_{t,\e}^c \giv A\ui_\e]   \\
           &\Le 2\e + \pr_{x_0}[A_{t,\e}^c \giv A\ui_\e].
\end{align*}
Hence, and from \Ref{ADB-CLT-8} and~\Ref{ADB-A1-prob},
\eq\label{ADB-main-approx-CLT}
    |\pr_{x_0}[\{\log X_t - \bb t \ge y\sqrt t\} \cap \cX] - (1-\Phi(y/v))\pr_{x_0}[\cX]| 
           \Le  2\pr_{x_0}[A_{t,\e}^c \giv A\ui_\e] + \d(t,\e) + 4\e.
\en

To bound $\pr_{x_0}[A_{t,\e}^c \giv A\ui_\e]$, note that, from~\Ref{ADB-At-def}, 
\[
     \pr_{x_0}[A_{t,\e}^c \giv A\ui_\e] \Le \pr_{x_0}[(E_\e\ui)^c \giv A\ui_\e] + \pr_{x_0}[(E_\e\ut)^c \giv A\ui_\e] 
                 + \pr_{x_0}[\tZ^*_t > \tz(\e) \giv A\ui_\e].
\]
From~\Ref{ADB-E1-def}, we have $\pr_{x_0}[(E_\e\ui)^c \giv A\ui_\e] \le 2\e$, and
$\pr_{x_0}[\tZ^*_t > \tz(\e) \giv A\ui_\e] < \e$, from~\Ref{ADB-tZs-prob}, if $t > 2t^* + s(\e)$.  It remains to bound 
$\pr_{x_0}[(E_\e\ut)^c \giv A\ui_\e]$.  From~\Ref{ADB-E2-def}, we have
\eqs
   \pr_{x_0}[(E_\e\ut)^c \giv A\ui_\e] &\le& \pr_{x_0}[(H(j_1(\e),\t(j_1(\e))) 
                                \cap  \cE^{\h}(j_1(\e),\t(j_1(\e))))^c \giv A\ui_\e]  \\
           &&\qquad\mbox{} +  \pr_{x_0}[\{\hcE^{\h}(j_1(\e),\t(j_1(\e)))\}^c \giv A\ui_\e], 
\ens
and, from Lemma~\ref{ADB-prob-HE}, we have
\[
    \pr_{x_0}\bigl[\bigl\{H(j_1(\e),\t(j_1(\e))) \cap  \cE^{\h}\bigl(j_1(\e),\t(j_1(\e))\bigr)\bigr\}^c \giv A\ui_\e\bigr]
            \Le C_0 j_1^{-\d_0}.
\]
From the definition~\Ref{ADB-E-hat-def} of~$\hcE^{\h}(j_1,\t(j_1))$, 
\[
    \pr_{x_0}[\{\hcE^{\h}(j_1(\e),\t(j_1(\e)))\}^c \giv A\ui_\e] 
              \Eq \pr\bigl[(\a/\bet')\inf_{i\ge1}(V_i - i\f') > -\half\log j_1(\e)\bigr],
\]
where~$(V_i,\,i\ge1)$ are as in Lemma~\ref{ADB-ex-min}, and where
\[
      \f' \Def \bet'(\log(1/p) + 3\h/2)/\a \ <\ 1.
\]
By an argument as for Lemma~\ref{ADB-ex-min}~(i), it follows that
\[
    \pr[\{\hcE^{\h}(j_1(\e),\t(j_1(\e)))\}^c \giv A\ui_\e] \Le \exp\{-(\half\log j_1(\e))\bet'u(\f')/\a\} 
              \Eq (j_1(\e))^{-\d_1},
\]
where $\d_1 > 0$.  Combining these bounds, we conclude that
\eq\label{ADB-E2-prob}
    \pr_{x_0}[(E_\e\ut)^c \giv A\ui_\e] \Le C_1 (j_1(\e))^{-(\d_0\wedge\d_1)}.
\en
Hence 
\[
    \pr_{x_0}[(A_{t,\e})^c \giv A\ui_\e] \Le 3\e + C_1 (j_1(\e))^{-(\d_0\wedge\d_1)}  \qquad\mbox{if}\quad t > 2t^* + s(\e).
\]
Thus, and from \Ref{ADB-main-approx-CLT} and~\Ref{ADB-CLT-6}, it follows that
\eq\label{ADB-last-line}
    \limsup_{t\to\infty} |\pr_{x_0}[\{\log X_t - \bb t \ge y\sqrt t\} \cap \cX] 
                   - (1-\Phi(y/v))\pr_{x_0}[\cX]|
\en
can be made arbitrarily small by choice of~$\e$, and is therefore equal to zero. 
\adbb{ 
Finally, if~$Y_t$ is any stochastic process such that $Y_t \to \infty$ a.s., then
$\lti\pr[Y_t \le M] = 0$ for all $M \in \re$.  Applying this observation using $\pr_{x_0}[\cdot \giv \cX^c]$
as probability measure
shows that, since $\log X_t \to -\infty$ on $\cX^c$, 
\[
   \limsup_{t\to\infty}\pr_{x_0}[\{\log X_t \ge  \bb t + y\sqrt t\} \cap \cX^c] 
            \Le \lti\pr_{x_0}[\log X_t \ge M \giv \cX^c] \Eq 0
\]
for any $y,M \in \re$,  completing the proof of the theorem.
} \ep
\end{proof}

\section{The degree weighted distribution}\label{ADB-degree-weighted}
\setcounter{equation}{0}

Throughout this section, we work with the process~$X$ under the restriction
\begin{align*}
     \alpha_k \Eq 1 - (1 - p_k)\beta_k  \Eq \a,\quad \text{for all $k\geq 1$,}
\end{align*}
for some fixed $\a > 0$,
which includes the generalized DD model of~\Ref{ADB-Q-matrix}. If $\mathbf{e}$ is defined to be the vector with 
$e_k=k$ for $k\geq 0$, then $Q\mathbf{e}=(2\alpha-1)\mathbf{e}$, implying that $\mathbf{e}$ is the 
$(1-2\alpha)$-invariant vector for~$Q$ on~$\mathbbm{N}$. 
\adbb{
Because of this, the matrix~$\tQ$ on~$\nat$ having elements
\begin{align}
       &\tQ_{k,k+1} \Eq \alpha(k+1);\qquad \tQ_{k,k} \Eq -\{-1+k\alpha + \beta_k(1-\pi_{kk})+2\alpha\}; \non\\
       &\tQ_{k,j} \Eq p_k\beta_k\frac{j\pi_{kj}}{kp_k}, \quad 1\leq j\leq k-1; \qquad \tQ_{k,j} \Eq 0,
                \quad j \ge k+2, \label{wtQmat}
 \end{align}
for each $k \ge 2$, and with $\tQ_{1,2} \Eq - \tQ_{1,1} = 2\a$, is a $Q$-matrix; indeed,
\begin{equation}\label{anotherexp}
    \tQ_{i,j}=(Q_{i,j}+(1-2\alpha)\delta_{ij})j/i\quad\text{for $i,j\in\mathbbm{N}$.}
\end{equation}
Note that, from the definition of~$p_k$, $\pi_{11} = p_1$, and $1 - (1-p_1)\bet_1 = \a$, so that the general formula
for~$\tQ_{k,k}$ in~\Ref{wtQmat}, when evaluated for $k=1$, still gives the correct value~$-2\a$.
The form of the matrix~$\tQ$ shows that the pure jump Markov process $\tX := (\tX_t,t\in \mathbbm{R}_+)$ on~$\nat$
corresponding to~$\tQ$ is dominated by a Yule process with rate~$2\a$, and is thus non-explosive.
Hence we can deduce the following proposition;
as observed in \cite{jordan}, in the case where $\alpha=p_k=p$ and $\Pi_k = \Bi(k,p)$, this follows from 
\textcite[Lemma 3.3]{pollett1988reversibility} when $2\a < 1$. 
} 

\begin{proposition}\label{quasi}
\adbb{
The point probabilities for the processes $X$ and~$\tX$ are related by
 \begin{equation}\label{2procs}
     j\pr_i[X_t=j] \Eq e^{-(1-2\alpha)t}\,i\, \pr_i[\tX_t=j]\quad\text{for all $i,j\in \mathbbm{N}$ and $t\geq 0$.}
 \end{equation}
}
\end{proposition}

\begin{proof}
\adbb{
Because~$X$ is a non-explosive pure jump Markov process, the Kolmogorov forward equations
\[
     \frac d{dt} \pr_i[X_t=j] \Eq  \sum_{l \ge 0}\pr_i[X_t=l]  Q_{l,j}
\]
are satisfied.   Hence, writing $x\ui_t(j) := e^{(1-2\a)t}\,(j/i)\,\pr_i[X_t=j]$, it follows that
\[
    \frac d{dt} x\ui_t(j) \Eq \sum_{l \ge 0} x\ui_t(l) (j/l) Q_{l,j} + (1-2\a) x\ui_t(j)
             \Eq \sum_{l \ge 0} x\ui_t(l) \tQ_{l,j},
\]
this last from~\Ref{anotherexp}.  Hence~$x\ui_t$ satisfies the forward equations for the process~$\tX$.
On the other hand, by \cite[Theorem~2]{hamza1995conditions}, the process 
\[
     \be(X_t) - \int_0^t (Q\be)(X_u)\,du \Eq X_t - \int_0^t (2\a-1) X_u\,du
\]
is a martingale, implying that $m\ui(t) := \ex\{X_t \giv X_0=i\}$ satisfies
\[
    m\ui(t) \Eq i + (2\a-1) \int_0^t m\ui(u)\,du,\qquad \mbox{or}\qquad m\ui(t) \Eq i e^{(2\a-1)t}.
\]
Hence, for all $t \ge 0$, 
\[
     \sum_{l \ge 0} x\ui_t(l) \Eq e^{(1-2\a)t}\sum_{l \ge 0} (l/i) \pr_i[X_t=l] \Eq e^{(1-2\a)t}i^{-1}m\ui(t) \Eq 1.
\]
Thus~$x\ui_t$ is a probability distribution for each~$t$, and hence is the  unique probability
solution to the forward equations for~$\tX$ having initial condition the point mass on~$i$;  that is,
\[
        x\ui_t(j) \Eq e^{(1-2\a)t}\,(j/i)\,\pr_i[X_t=j] \Eq \pr_i[\tX_t = j],
\]
proving the proposition. \ep
}
\end{proof}

\adbb{
The matrix~$\tQ$ has the same form as that given in~\Ref{ADB-Q-matrix-1}, with $\tal_k := \a(k+1)/k$ and
$\tbe_k := p_k\bet_k$, and with $\tPi_k$ defined by $\tpi_{kj} := j\pi_{kj}/kp_k$, and in particular with $\tpi_{11} = 1$.
Note that the quantities~$\tal_k$ are not the same for each~$k$, but that differing~$\a_k$ were
allowed for in~\Ref{ADB-Q-matrix-1}.  
The only effect of not having $\tpi_{11} < 1$ is to make~$\nat$ a closed class, and to remove the absorbing state at~$0$.
The mean and variance of~$\tPi_k$ are given by
\[
     k\tp_k \Def \sum_{j=1}^{k} j\tpi_{kj} \Eq \sum_{j=1}^{k} \frac{j^2\pi_{kj}}{kp_k}
         \andy \tsi_k \Def \sum_{j=1}^{k} \frac{j^3\pi_{kj}}{kp_k} - (k\tp_k)^2.
\]
If the assumptions in \Ref{ADB-parameter-assns-again} and~\Ref{ADB-parameter-assns-gen} are satisfied for~$X$
(the condition on~$\a_k$ is automatic, since in this section $\a_k = \a$ for all~$k$),
then
\eq\label{ADB-tilde-assns-1}
    |\tbe_k - p\bet| \Le \tc_1 k^{-\g} \andy |\tal_k - \a| \Le \tc_4 k^{-1},
\en
so that the assumptions~\Ref{ADB-parameter-assns-gen} are satisfied fo~$\tX$ with $\tbe := p\bet$ and $\tal = \a$.
That assumptions~\Ref{ADB-parameter-assns-again} are satisfied fo~$\tX$, with $\g/4$ for~$\g$, is established
in the following lemma. 
} 

\begin{lemma}\label{sbvar}
Let $\tY_k \sim \tPi_k$. Then there are positive constants $\tc_2$ and~$\tc_3$
\begin{equation}\label{newvar}
   |\tp_k - p| \Le  \tc_2 k^{-\g}\quad \text{and} \quad k^{-2}\tsi^2_k \Le \tc_3 k^{-\gamma/4},
\end{equation}
\end{lemma}

\begin{proof}
It follows by direct calculation that 
\[
    k\tp_k \Eq \E\tY_k \Eq  kp_k + \sigma^2_k/(kp_k),
\]
from which, and from~\Ref{ADB-parameter-assns-again}, the first assertion follows. To prove the second assertion, it is
enough to show that $(kp_k)^{-2}\E \tY^2_k\leq 1 + O(k^{-\gamma/4})$, because 
$\tsi^2_k = \E(\tY^2_k) - (\E \tY_k)^2$, and $\E \tY_k=kp_k + O(k^{1-\g})$, from~\Ref{ADB-parameter-assns-again}. 
Letting $\varepsilon_k:=k^{-\gamma/4}$ and $Y_k \sim \Pi_k$, it follows from Chebyshev's inequality that
\begin{align*}
    (kp_k)^{-2}\E \tY^2_k &\Eq (kp_k)^{-3} \sum^k_{j=1}j^3\pi_{kj}\\
    &\Le (kp_k)^{-3}\bigl\{(1+\varepsilon_k)^3(kp_k)^3 + k^3\Prob[Y_k>(1+\varepsilon_k) kp_k]\bigr\}\\
    &\Le (1+\varepsilon_k)^3 + O(k^{-\gamma/2})\\
    &=1 + O(k^{-\gamma/4}),
\end{align*}
which completes the proof.
\end{proof}

\adbb{
As a result of the inequalities~\Ref{ADB-tilde-assns-1} and of Lemma~\ref{sbvar}, we can now apply the
theorems of the previous sections to the process~$\tX$.  Recall that, as in~Theorem~\ref{ADB-recurrence-2},
for any $u > 1$, $x^*(u) < u$ denotes the positive solution to the equation $x = u(1-e^{-x})$.
}

\begin{theorem}\label{recurrence} 
For the degree weighted continuous time process~$\tX$, we have the following asymptotic behaviour:
\begin{enumerate}
\item  If $\alpha<p\beta\log(1/p)$, then $\tX$ is geometrically ergodic, and the stationary distribution 
has $\eta$-th moment finite for all $\h < \h^* := x^*(p\b\log(1/p)/\a)/\log(1/p) < p\b/\a$, and infinite
for all $\h > \h^*$;
\item  If $\alpha=p\beta\log(1/p)$, then $\tX$ is null recurrent; 
\item  if $\alpha>p\beta\log(1/p)$, then $\tX$ is transient. 
\end{enumerate}
\end{theorem}

\adbb{
In case~1,  $\pr_i[\tX = j] \to \tp_j(\infty)$ for each $j \ge 1$, where~$\tp(\infty)$ is a probability
distribution on~$\nat$.  In view of Proposition~\ref{quasi}, it follows that
\[
 e^{(1-2\alpha)t}\pr_i[X_t=j] \Eq (i/j) \pr_i[\tX_t=j] \ \to\ (i/j) \tp_j(\infty) \quad \mbox{as}\ t \to \infty,
\]
so that the probabilities~$\pr_i[X_t=j]$ all decay exponentially in time, with the same rate.
Furthermore, since the function $j \mapsto j^{-1}$ is bounded on~$\nat$,
\eq\label{ADB-exponential-decay}
   e^{(1-2\alpha)t}\pr_i[X_t \in A] \Eq i\sum_{j \in A} j^{-1}\pr_i[\tX_t=j] \ \to\ i\sum_{j \in A} j^{-1} \tp_j(\infty),
\en
and hence
\eq\label{ADB-conditional-convergence}
    \pr_i[X_t = j \giv X_t \ge 1] \Eq \frac{j^{-1}\pr_i[\tX_t=j]}{\sum_{l \ge 1} l^{-1}\pr_i[\tX_t = l]}
             \ \to\ \frac{j^{-1}\tp_j(\infty)}{\sum_{l \ge 1} l^{-1} \tp_l(\infty)},
\en
showing that then the conditional distribution of~$X_t$, given $X_t \ge 1$, converges, and identifying the limit.
In the basic model, with $q=0$, \Ref{ADB-conditional-convergence} is \textcite[Proposition 3.4]{jordan}. 
and the exponential tail in (\ref{ADB-exponential-decay}) is in \textcite[Corollary 2.7(1)]{hermann2020markov}. 
}

\begin{figure}[ht]
\centering
\includegraphics[width=9cm]{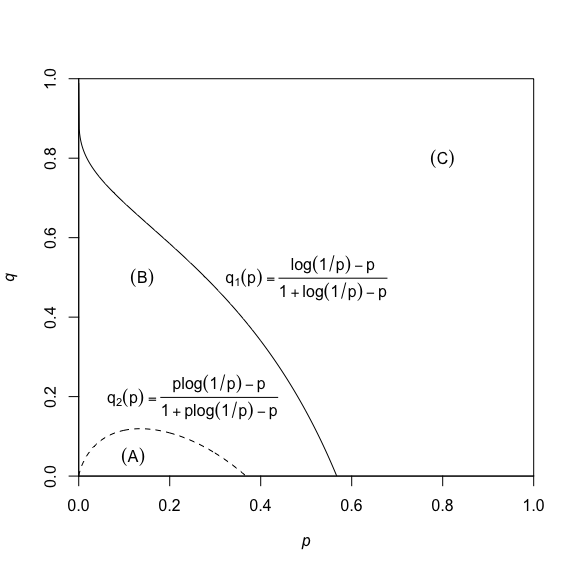}
\caption{\small A $p-q$ plane. The labelled areas correspond to the cases above. }
\label{TL-pq-picture}
\end{figure}

In the basic model, in which  $\alpha=p(1-q)+q$ and $\b = 1-q$, for any choices of $p$ and~$q$ such that $0<p<1$ and 
$0\leq q < 1$, the behaviour given in Theorems \ref{ADB-hX-recurrence} and~\ref{recurrence} can be categorized
according to the cases 
\eqs
    &&{\rm (A)\colon}\ \alpha<p(1-q)\log(1/p), \quad
      {\rm (B)\colon}\ p(1-q)\log(1/p)\leq \alpha <(1-q)\log(1/p) \\
    &&\qquad\andy {\rm (C)\colon}\ \alpha > (1-q)\log(1/p);
\ens 
see Figure~\ref{TL-pq-picture}.  For fixed~$p$, these cases
can be represented by (A): $q < q_2(p)$, (B): $q_2(p) < q < q_1(p)$ and (C): $q > q_1(p)$,
where~$q_1$ is the inverse of the function~$p_*$ solving~\Ref{ADB-p(q)-equation}:
\[
    q_1(p) \Eq \frac{\log(1/p) - p}{1 + \log(1/p) - p} \ >\ 
    q_2(p) \Eq \frac{p\log(1/p) - p}{1 + \log(1/p) - p},     
\]
and where $q_2(p) < 0$ if $p > e^{-1}$ and $q_1(p) < 0$ if $p > p_*(0) \approx 0.5671$.

In particular, $q_2(0) = q_2(e^{-1}) = 0$, and the maximum of~$q_2$ is given by $q_2(e^{-2}) = 1/(e^2 + 1) \approx 0.1192$. 
Case~(A) holds for $0<p<e^{-1}$ and $0<q<q_2(p)$, and in this case \Ref{ADB-exponential-decay} 
and~\Ref{ADB-conditional-convergence} hold. In case~(B), the process~$X$ is absorbed in zero with
probability one, but the degree weighted process is transient.  In case~(C), the process~$X$ is transient
and the central limit theorem of Theorem~\ref{CLT} can be used to approximate its (time-dependent)
distribution.

When the process~$\tX$ is transient, Theorem \ref{aslimit} and~\ref{CLT} can be applied to it,
in view of the inequalities~\Ref{ADB-tilde-assns-1} and of Lemma~\ref{sbvar}.  However, because of
the degree weighting, their conclusions no longer give information about the {\it typical\/} degree
distribution.

\section{Discrete time processes}\label{ADB-discrete-time}
\setcounter{equation}{0}

\subsection{The basic DD model}
We now return to the analysis of the original probabilities~$\pp_{\snt}\TT$, $\snt \in \Z_+$,
that can be derived as the state probabilities
of an inhomogeneous Markov chain~$Y := (Y_{\snt},\,\snt \in \Z_+)$, using~\Ref{ADB-Q-version-discrete}.   
Note first that, given that the chain~$Y$ has a jump at step $(\snt+1)$, the probability distribution
of~$Y_{\snt+1}$, given~$Y_{\snt} = k$, depends only on the value~$k$, and is the same for all~$s$; the distribution is given by
the transition probabilities out of the state~$k$ of the {\it jump chain\/} associated with
the homogeneous process~$X$ in continuous time.  Hence the Markov chain~$Y$ can be constructed from
a realization of the jump chain of~$X$, together with an independent realization of a sequence $(U_j,\,j\ge1)$
of independent standard uniform random variables, which are used to determine the residence times in the
successive states of the chain.  \adbb{In particular, properties such as recurrence and transience, which 
depend only on the sequence of states visited, can be deduced for~$Y$ from the corresponding properties for~$X$.}
However, the factors $(\snt+1)^{-1}$
in the definition of the jump probabilities at time~$\snt$ imply that the distribution of the residence times 
of~$Y$ in a given state~$k$ themselves depend on the time~$\snt$ at which~$k$ was reached, and that the residence times,
as measured by the number of steps in~$k$, become much longer as~$\snt$ increases.  
Instead, by compressing the time
scale, a close analogue~$\Yh$ of the homogeneous process~$X$ can be derived, as follows.

Let~$(\hX_j,\,j\ge0)$ be a realization of the jump chain of~$X$.  Then a version of~$X$ can be obtained 
from~$\hX$ and an independent sequence
of mutually independent standard uniform random variables $(U_j,\,j\ge1)$ by setting
\eq\label{ADB-X-construct}
      X_t \Eq \hX_n \ \quad\mbox{if}\ \quad 
             S_n \Def \sum_{j=1}^{n} V_j \Le t \ <\ S_{n+1},\quad t \ge 0,
\en
where, with $\hq_k := -Q_{k,k}$, $V_j := \hq_{\hX_{j-1}}^{-1}\{-\log(1-U_j)\}$ represents the $j$-th residence time of~$X$, 
the residence time following its $(j-1)$-st jump, and has an exponential distribution with mean $\hq_{\hX_{j-1}}^{-1}$. 
\adbb{Note that, if $\hX_{j-1} = 0$, then $V_j = \infty$ a.s., because $q_0=0$.}  
A similar construction can be used to define~$\Yh$.
Let~$N_j$ represent the step at which the $j$-th jump of~$Y$ occurs, with $N_0 = \snt_0$. 
 For any $a > -1$, define the \adbb{generalized harmonic numbers $h_a(j)$, $j\ge1$, by 
$h_a(j) := \sum_{l=1}^j 1/(l+a)$, so that~$h_0$ represents the harmonic numbers~$h$ defined in~\Ref{ADB-harmonic-numbers}.
Let} $\tV_j := h_0(N_j) - h_0(N_{j-1})$ represent the $j$-th residence time of~$Y$, 
but in an almost logarithmically distorted time scale; \adbb{again, if $\hX_{j-1} = 0$, then $\tV_j = \infty$ a.s.}  Then set
\eq\label{ADB-tY-construct}
    \Yh_t \Def \hX_n  \ \quad\mbox{if}\ \quad  \tS_n \Def \sum_{j=1}^{n} \tV_j \Le t \ <\ \tS_{n+1},
               \quad t \ge 0.
\en
From \Ref{ADB-X-construct} and~\Ref{ADB-tY-construct}, it is immediate that the paths of $X$ and~$\Yh$ are
close to one another if the partial sum processes $S$ and~$\tS$ defined from the $V_j$ and the~$\tV_j$ are close to one
another.
\adbb{
Indeed, defining
\[
      \D_0 \Def 0;\qquad \D_n \Def S_n - \tS_n, \quad n \ge 1,
\]
and then setting
\eq\label{ADB-Delta-def} 
     \D(t) \Def \frac{\tS_{n+1} - t}{\tS_{n+1} - \tS_n}\,\D_n + \frac{t - \tS_n}{\tS_{n+1} - \tS_n}\,\D_{n+1} 
               \ \quad\mbox{if}\ \quad \tS_n \Le t \ <\ \tS_{n+1},\qquad t \ge 0,
\en
\adbb{with $\D(t) := \D_n$ if $t \ge \tS_n$ and $\tS_{n+1} = \infty$},
we have
\eq\label{ADB-matching-processes}
   \Yh_t \Eq X_{t + \D(t)}, \quad t \ge 0, 
          \qquad\mbox{and also}\qquad Y_{\snt} \Eq \Yh_{h_0(\snt) - h_0(\snt_0)}, \quad \snt \ge \snt_0 \in \Z_+.
\en
In this way, the process~$\Yh$ is represented as a random time shift of the process~$X$, and 
the process~$Y$ can be deduced from~$\Yh$ by a deterministic time change.  We now show that
the time shift connecting $\Yh$ and~$X$ is not substantial:  in fact, 
\eq\label{ADB-Delta-limit}
    \D \Def \lim_{t\to\infty}\D(t) \quad \mbox{exists a.s.}
\en
}

To establish~\Ref{ADB-Delta-limit}, we need to examine the distribution of the random variables~$\tV_j$.
From the definition of the process~$Y$, conditional on the values of the whole of~$\hX$ and on $\{N_l,\,1 \le l\le j-1\}$, 
and on the event $\{(N_{j-1},\hX_{j-1})=(\snt,k)\}$,
$\tV_j$ has the distribution of $V_0(\snt,\hq_{k})$, where  $V_a(\snt,\bsb)$ denotes a random variable taking values 
in the discrete set $\{h_a(j) - h_a(\snt),\ j \ge \snt+1\}$ with probabilities given by
\eq\label{ADB-Snk-def}
   \pr[V_a(\snt,\bsb) > h_a(r) - h_a(\snt)] \Eq \prod_{l=\snt}^{r-1} (1 - \bsb/(l+a+1)),\qquad r \ge \snt.
\en
The next lemma shows that if~$E_\bsb$ denotes a random variable having an exponential distribution with mean~$1/\bsb$, 
random variables with the distributions of $V_a(\snt,\bsb)$ and~$E_\bsb$ can be 
constructed on the same probability space, in such a way that they are close.

\begin{lemma}\label{ADB-V-E-coupling}
 For $V_a(\snt,\bsb)$ defined as above, $\ex V_a(\snt,\bsb) = 1/\bsb$.
Furthermore, if~$E_\bsb$, as defined above, is coupled to~$V_a(\snt,\bsb)$ 
using the quantile coupling, and if $\bsb/(\snt+a+1) \le \f < 1$, then there exists a positive constant~$c_\f$
such that, for all $j \ge \snt+1$,
\eq\label{ADB-S-E-ineq}
    -\frac1{j+a} \Le E_\bsb - V_a(\snt,\bsb) \Le c_\f\,\frac{\bsb(h_a(j) - h_a(\snt))}{\snt+a+1} \,,
\en
almost surely on the event $\{V_a(\snt,\bsb) = h_a(j) - h_a(\snt)\}$.  As a consequence, there exists a positive 
constant~$k_{\f}$ such that
\[
    \ex\{(E_\bsb - V_a(\snt,\bsb))^2\} \Le \frac {k_{\f}}{(\snt+a+1)^{2}}\,.
\]
\end{lemma}

\begin{remarka}\label{ADB-quantile-coupling}
{\rm The quantile coupling can be achieved by using a standard uniform random variable~$U$:  
take $E_\bsb := -\bsb^{-1}\log(1-U)$ and}
\[
    V_a(\snt,\bsb) \Eq h_a(r+1) - h_a(\snt) \quad\mbox{if}\quad 
           \prod_{l=\snt}^{r} (1 - \bsb/(l+a+1))  \Le 1-U \ <\ \prod_{l=\snt}^{r-1} (1 - \bsb/(l+a+1)), \quad 
\]
for $r \ge \snt$.
\end{remarka}

\begin{proof}
\adbb{
For any non-negative non-decreasing sequence $(f_j,\, j\ge1)$ and any probability distribution $(p_j,\,j\ge1)$
on~$\nat$, it is immediate that, if $\sum_{j\ge1}p_j f_j < \infty$, then
\[
      \sum_{j\ge1}p_j f_j \Eq \sum_{j \ge 1} P_j(f_j - f_{j-1}),        
\]
where $P_j := \sum_{l \ge j}p_l$, and $f_0 := 0$.  Hence, 
}
since $h_a(j+1) - h_a(j) = 1/(j+a+1)$, it follows that
\eqs
    \ex V_a(\snt,\bsb) &=& \sum_{j \ge \snt} \frac1{j+a+1}\,\prod_{l=\snt}^{j-1} (1 - \bsb/(l+a+1)) \\
       &=& \frac1\bsb\, \sum_{j \ge \snt} 
                  \Bigl\{ \prod_{l=\snt}^{j-1} (1 - \bsb/(l+a+1)) - \prod_{l=\snt}^{j} (1 - \bsb/(l+a+1)) \Bigr\} \\
       &=& \frac1\bsb\,.
\ens

For the quantile coupling, first observe that
\[
    \prod_{l=\snt}^{j-1} (1 - \bsb/(l+a+1)) \Eq \exp\Bigl\{ \sum_{l=\snt}^{j-1} \log(1 - \bsb/(l+a+1)) \Bigr\},
\]
from which it follows that
\eqa
   e^{-\bsb(h_a(j) - h_a(\snt))} &\ge& \pr[V_a(\snt,\bsb) > h_a(j) - h_a(\snt)] \label{ADB-S-E-1}\\ 
     &=&  e^{-\bsb(h_a(j) - h_a(\snt))}\,
        \exp\biggl\{ \sum_{l=\snt}^{j-1} \Bigl\{\log\Bigl(1 - \frac{\bsb}{l+a+1}\Bigr) + \frac{\bsb}{l+a+1}\Bigr\} \biggr\} \non\\
     &=& e^{-\bsb(h_a(j) - h_a(\snt) + \ps_a(j,\bsb))} \Eq \pr[E_b > h_a(j) - h_a(\snt) + \ps_a(j,\bsb)],
               \label{ADB-S-E-2}  
\ena
say.
Because $\bsb \le \f(\snt+a+1)$ \adbb{and $\log(1-x) > -x - c_\f x^2$ in $0 < x < \f$, for some constant~$c_\f > 0$,} 
it is immediate that
\eqs
   \ps_a(j,\bsb) &=& 
            -\frac1\bsb\,\sum_{l=\snt}^{j-1} \Bigl\{\log\Bigl(1 - \frac{\bsb}{l+a+1}\Bigr) + \frac{\bsb}{l+a+1}\Bigr\} \\
            &\le& \frac{c_\f}{\bsb}\, \sum_{l=\snt}^{j-1} \frac{\bsb^2}{(l+a+1)^2} 
            \Le c_\f\,\frac{\bsb(h_a(j) - h_a(\snt))}{\snt+a+1}\,.
\ens
Hence, using the quantile coupling, if $V_a(\snt,\bsb) \le h_a(j) - h_a(\snt)$, it follows that
\[
    E_\bsb \Le h_a(j) - h_a(\snt) + \ps_a(j,\bsb) \Le h_a(j) - h_a(\snt) + c_\f\,\frac{\bsb(h_a(j) - h_a(\snt))}{\snt+a+1}\,,
\]
and, in particular, on $\{V_a(\snt,\bsb) \le h_a(j) - h_a(\snt)\}$, 
the second inequality in~\Ref{ADB-S-E-ineq} now follows.
The first inequality in~\Ref{ADB-S-E-ineq} follows similarly as a consequence of~\Ref{ADB-S-E-1}, applied with $j-1$ for~$j$.

For the final statement, it follows from \Ref{ADB-S-E-ineq} and~\Ref{ADB-S-E-1} that
\eqs
     \lefteqn{\ex\{(E_\bsb - V_a(\snt,\bsb))^2\}} \\
    && \Le
         \sum_{j \ge \snt+1} \pr[V_a(\snt,\bsb) = h_a(j) - h_a(\snt)] 
               \max\biggl\{ \frac1{(\snt+a+1)^2},  \Bigl(c_\f\,\frac{\bsb(h_a(j) - h_a(\snt))}{\snt+a+1}\Bigr)^2 \biggr\} \\
     && \Le  \frac1{(\snt+a+1)^2} 
         + \sum_{j \ge \snt+1} \frac{\bsb}{j+a}\,  e^{-\bsb(h_a(j-1) - h_a(\snt))} \Bigl(c_\f\,\frac{\bsb(h_a(j) - h_a(\snt))}{\snt+a+1}\Bigr)^2 \\
     && \Le  \frac1{(\snt+a+1)^2} 
         + \sum_{j \ge \snt} \frac{2e\bsb}{j+a+1}\,  e^{-\bsb(h_a(j) - h_a(\snt))} \Bigl(c_\f\,\frac{\bsb(h_a(j) - h_a(\snt))}{\snt+a+1}\Bigr)^2.
\ens
Now, for a unimodal function $f\colon \re_+ \to \re_+$ and for a sequence $0 = x_0 < x_1 < \cdots$, with
$x_{j+1} - x_j \le x_j - x_{j-1}$ for all $j \ge 1$, it follows by an elementary comparison that
\[
     \sum_{j \ge 0}(x_{j+1} - x_j) f(x_j) \Le \int_0^\infty f(x)\,dx + 2(x_1 - x_0)\sup_{x \ge 0}f(x).
\]
Hence
\eqs
 \lefteqn{ \sum_{j \ge \snt} \frac{2e\bsb}{j+a+1}\,  e^{-\bsb(h_a(j) - h_a(\snt))} 
                            \Bigl(c_\f\,\frac{\bsb(h_a(j) - h_a(\snt))}{\snt+a+1}\Bigr)^2 }\\
     &&\Le \frac{2e c_\f^2}{(\snt+a+1)^2}\, \Bigl\{ \int_0^\infty x^2 e^{-x}\,dx + 8\adbb{\f} e^{-2} \Bigr\},
\ens
and the final part of the lemma follows.\ep
\end{proof}

This enables the partial sum processes defined from the $V_j$ and the~$\tV_j$ to be constructed so as to be close to one
another.

\begin{corollary}\label{ADB-path-coupling}
Conditional on the whole of the jump chain~$\hX$, the random variables $(\tV_j,\,j\ge1)$ and~$(V_j,\,j\ge1)$ can be
coupled in such a way that
\adbb{
\[
      \lnti \D_n \Eq \lti \D(t) \Eq \D
\]
exists a.s., with the convergence holding also in mean square.}
\end{corollary}

\begin{proof}
In view of the constructions \Ref{ADB-X-construct} and~\Ref{ADB-tY-construct}, it is enough to use the coupling
in Remark~\ref{ADB-quantile-coupling} \adbb{with $a=0$} successively, with a sequence of independent standard uniform random
variables $(U_j,\,j\ge1)$, that are also independent of the jump chain~$\hX$.
\adbb{
The corollary then follows because $(\D_r,\,r\ge0)$ is a square integrable martingale,
as is implied by the variance bound given in Lemma~\ref{ADB-V-E-coupling}, and because
$\D(\cdot)$ interpolates the values of $(\D_r,\,r\ge0)$.
}
Note that, under the assumptions of Section~\ref{intro}, $\hq_k \le (\snt+1)\f$ for all
$k \le \snt$, where $\f := \half(1 + \max_{k\ge1}\a_k) < 1$,
and that necessarily $\hX_j \le N_j$ for all~$j$, because upward jumps are of magnitude~$1$. \ep
 \end{proof}

\adbb{
As a result of Corollary~\ref{ADB-path-coupling}, it is possible to transfer results about the
long term behaviour of~$X$ into corresponding results about~$\Yh$, and hence about~$Y$.  
When~$X$ is absorbed in zero with probability one, as in Corollary~\ref{ADB-recurrence}, the same
is true for $\Yh$ and~$Y$, but this is already clear, because they have the same jump chain as~$X$.
If the process~$X$ is transient, then there is an almost sure growth theorem for~$Y$ that
can be deduced from Theorem~\ref{aslimit}.  To state it, let~$J_{\snt}$ denote the number of `downward'
jumps of~$Y$ up to time $\snt \ge \snt_0$, and let $\Jh_t$ denote the number of downward jumps of~$\Yh$
up to time~$t$, where a downward jump is a jump that increases~$Z$ by one in the process~$(X,Z)$;
since the jump chains of the three processes $X, \Yh$ and~$Y$ are identical, their downward jumps can also
be defined to be the same.  Note that $J_{\snt} = \Jh_{h_0(\snt) - h_0(\snt_0)}$ for all $\snt \ge \snt_0$.
}

\begin{theorem}\label{aslimit-disc}
\adbb{If $\a > \bet\log(1/p)$, let $W := \lim_{t\to\infty} e^{-\alpha t}p^{-Z_t} X_t$ be the almost sure limit
given in Theorem~\ref{aslimit}.  Then $\snt^{-\a} p^{-J_{\snt}} Y_{\snt}$ converges a.s.\ as $\snt \to \infty$ to a limit~$W'$
such that
\[
    W' \ =_d\ e^{\a\{\D  + (\log \snt_0 - h_0(\snt_0) + \g)\}} W ,                                  
\]
where~$\D := \lim_{t\to\infty} \D(t)$, and~$\g$ is Euler's constant.  Note that $W$ and~$\D$ are dependent
random variables, and that $\lim_{\snt_0\to\infty}(\log \snt_0 - h_0(\snt_0) + \g) = 0$.}
\end{theorem}

\begin{proof}
\adbb{ From the definition of~$W$,  it is immediate that
\[
    e^{-\a(t+\D(t))} p^{-Z_{t+\D(t)}} X_{t+\D(t)} \ \to\ W \quad \mbox{a.s. as } t \to \infty.
\]
Then, because $\Yh$ and~$X$ can be constructed together in such a way that $\Yh(t) = X_{t+\D(t)}$ for all $t \ge 0$,
it follows that
\[
     e^{-\a(t+\D(t))} p^{-\Jh_t} \Yh_t \ \to\ W \quad \mbox{a.s. as } t \to \infty.
\]
Taking $t = h_0(\snt) - h_0(\snt_0)$ for $\snt \ge \snt_0$, $\snt \in \nat$, and because $\D(t) \to \D$ a.s.\ as $t\to\infty$, 
this implies that
\[
     e^{-\a(h_0(\snt) - h_0(\snt_0))} p^{-J_{\snt}} Y_{\snt} \ \to\ e^{\a\D} W \quad \mbox{a.s. as } s \to \infty,
\]
and the theorem follows, because $\lim_{\snt\to\infty} \{h_0(\snt) - \log \snt\} = \g$.} \ep
\end{proof}

The logarithm of the inhomogeneous discrete time process~$Y$ also satisfies a central limit theorem.

\begin{theorem}\label{CLT-disc}
\adbb{ If $\bb := \a - \bet\log(1/p) > 0$, for any $x_0 \in \nat$, the process~$Y$ 
satisfies
\[
     \lim_{\snt\to\infty} \pr[\{\log \snt\}^{-1/2}\{\log Y_{\snt} - \bb \log \snt\} \ge y \giv Y_{\snt_0} = x_0] 
              \Eq  \pr_{x_0}[\cX](1 - \Phi(y/v)),
\]
where $v^2 := \bet\{\log(1/p)\}^2$ and $\pr_{x_0}[\cX] := \pr[\lti X_t = \infty \giv X_0 = x_0]$.
}
\end{theorem}

\begin{proof}
\adbb{
 We begin by proving the central limit theorem for~$\Yh$, arguing using the coupling $\Yh_t = X_{t + \D(t)}$,
$t \ge 0$, and deducing the theorem from Theorem~\ref{CLT}.
}

\adbb{
Writing $\log X_t = \a t - Z_t \log(1/p) - \log W_t$, where $\log W_t \to \log W > -\infty$ a.s.\ on~$\cX$,
and $\log W_t \to -\infty$ a.s.\ on~$\cX^c$, we have
\eq\label{ADB-Tiffany-1}
    |\log X_{t+\D(t)} - \log X_t| \Le \a|\D(t)| + \log(1/p)|Z_{t+\D(t)} - Z_t| + |\log W_{t+\D(t)} - \log W_t|.
\en
For $s \ge 1$ and $t \ge 2(s+1)$, define the events
\eqa
   A_s\ui &:=& \{|\D(u) - \D| \le 1 \ \mbox{for all}\ u \ge s\};\qquad
   A_t\uth \Def \{|Z_{t + t^{1/8}} - Z_{t - t^{1/8}}| \le t^{1/4}\} ; \non\\
   A_t\ut &:=& \{|\D| \le t^{1/8}\}; \qquad   
   A_s\uf \Def \{|\log W_u - \log W_v| \le 1 \ \mbox{for all}\ u,v \ge s\}.  \label{ADB-Tiffany-1.5}  
\ena
The first three of these events are shown to have probabilities approaching~$1$ as $s,t \to \infty$, and the fourth event
approaches the event~$\cX$ as $s \to \infty$.  On $A_s\ui \cap A_t\ut \cap A_s\uf$, with $s \ge 1$ and $t \ge 2(s+1)$,
so that $t - t^{1/8} - 1 > s$, it follows that $|\log W_{t+\D(t)} - \log W_t|\le 1$.  Hence, if also~$A_t\uth$
holds, then, using~\Ref{ADB-Tiffany-1},
\[
    |\log X_{t+\D(t)} - \log X_t| \Le \a(t^{1/8} +1) + \log(1/p) t^{1/4} + 1 \ =:\ r(t);
\]
here, we have used the fact that~$Z_u$ is a.s.\ non-decreasing in~$u$.  Note that $t^{-1/2}r(t) \to 0$
as $t\to\infty$, so that, on the event $\cE_{s,t} := A_s\ui \cap A_t\ut \cap A_t\uth \cap A_s\uf$, the difference
between $\log X_{t+\D(t)}$ and~$\log X_t$ is negligible as far as the CLT is concerned.  We exploit this as follows.
}

\adbb{
First, since $\log X_t \to -\infty$ a.s.\ as $t\to\infty$ on~$\cX^c$, it follows also that $\log X_{t+\D(t)} \to -\infty$
a.s.\ on~$\cX^c$.  Hence, as in concluding the proof of Theorem~\ref{CLT},
\eq\label{ADB-Xc-prob}
   \lti \pr_{x_0}[\{t^{-1/2}(\log X_{t+\D(t)} - \bb t) \ge y\} \cap \cX^c] \Eq 0.
\en
\ignore{
we show that the probability $\pr_{x_0}[\{t^{-1/2}(\log X_{t+\D(t)} - \bb t) \ge y\} \cap \cX^c]$ tends to zero
as $t\to\infty$.  For any positive $a,b$ and~$h$, we have
\eqs
   &&\pr_{x_0}[\{X_{t+\D(t)} > a\} \cap \cX^c] \Le
     \pr_{x_0}\bigl[\bigl\{\sup_{t-h \le u \le t+h} X_u > a\bigr\} \cap \cX^c\bigr] + \pr_{x_0}[|\D(t)| > h] \\
  &&\quad\Le \pr_{x_0}\bigl[\bigl\{\sup_{t-h \le u \le t+h} X_u > a\bigr\} \cap \{X_{t-h} \le b\}\bigr] 
      + \pr_{x_0}[\{X_{t-h} > b\} \cap \cX^c] + \pr_{x_0}[|\D(t)| > h] .
\ens
For the first probability, it is enough to observe that the growth of~$X$ is stochastically bounded above by
that of a Yule process~$X^y$
with rate~$\a$, for which, being a non-decreasing process, $\sup_{t-h \le u \le t+h} X^y_u = X^y_{t+h}$.
Hence, since $\ex\{X^y(t+h) \giv X^y_{t-h} \le b\} \le b e^{2\a h}$, it follows by Markov's inequality that
\[
   \pr_{x_0}\bigl[\bigl\{\sup_{t-h \le u \le t+h} X_u > a\bigr\} \cap \{X_{t-h} \le b\}\bigr]  \Le \frac{be^{2\a h}}a\,.
\]
Then $\pr_{x_0}[\{X_{t-h} > b\} \cap \cX^c] \Le C_0 b^{-\d_0}$, by Lemma~\ref{ADB-prob-HE}.  Taking $a = e^{\bb t/2}$,
$b = \sqrt a$ and~$h = 1 + t^{1/8}$, it follows easily that, for any real~$y$,
\eq\label{ADB-Xc-prob}
  \lim_{t \to \infty}\pr_{x_0}[\{t^{-1/2}(\log X_{t+\D(t)} - \bb t) \ge y\} \cap \cX^c] 
       \Le \limsup_{t\to\infty}\{\pr_{x_0}[(A_t\ui)^c] +\pr_{x_0}[(A_t\ut)^c]\} \Eq 0,
\en
because $\D(t) \to \D$ a.s.\ and because~$\D$ is a fixed random variable.
}
It thus remains to approximate the probability $\pr_{x_0}[\{t^{-1/2}(\log X_{t+\D(t)} - \bb t) \ge y\} \cap \cX]$.
For this, we use the following sandwich:
\eqa
   \lefteqn{ \bigl(\{t^{-1/2}(\log X_t - \bb t) \ge y + t^{-1/2}r(t)\} \cap \cX\bigr) \setminus (\cE_{s,t}^c \cap \cX) }\non\\
    &&\ \subset\ \bigl(\{t^{-1/2}(\log X_{t+\D(t)} - \bb t) \ge y\} \cap \cX\bigr) \label{ADB-Tiffany-2}\\
    &&\ \subset\ \bigl(\{t^{-1/2}(\log X_t - \bb t) \ge y - t^{-1/2}r(t)\} \cap \cX\bigr) \cup (\cE_{s,t}^c \cap \cX). \non
\ena
}

\adbb{
The first step now is to show that the event
\[
   \cE_{s,t}^c \cap \cX \ \subset\ (A_s\ui)^c \cup (A_t\ut)^c 
            \cup (A_t\uth)^c \cup (\cX \setminus A_s\uf) 
\]
has small probability when $s$ and~$t$ are large.  It is immediate that 
\[
     \lim_{s \to \infty}\pr_{x_0}[(A_s\ui)^c] \Eq \lim_{t \to \infty}\pr_{x_0}[(A_t\ut)^c] \Eq 0,
\]
because $\D(t) \to \D$ a.s.\ and because~$\D$ is a.s.\ finite.
If $\b^* := \sup_{k \ge 0}\b_k$, a comparison between~$Z$ and a Poisson process with rate~$\b^*$ shows that
\[
    \pr_{x_0}[(A_t\uth)^c] \Le \Po(\b^*(t^{1/8}+1))\{(t^{1/4},\infty)\} \ \to\ 0 \quad \mbox{as } t \to \infty.
\]
For $\cX \setminus A_s\uf$, note that $A_s\uf$ is an increasing sequence, and that $\cX \subset \lim_{s\to\infty}A_s\uf$,
because $\cX = \{\lti \log W_t \in \re\}$.
On the other hand, for all $s \ge 1$,
\[
     \cX^c \Eq \{\lti \log W_t = -\infty\} \ \subset\ (A_s\uf)^c.
\]
Hence $\cX = \lim_{s\to\infty}A_s\uf$ a.s., and so $\lim_{s\to\infty} \pr_{x_0}[\cX \setminus A_s\uf] = 0$.
Thus it follows that 
\eq\label{ADB-limlim}
     \lim_{s \to \infty} \lti \pr_{x_0}[\cE_{s,t}^c \cap \cX] \Eq 0.
\en
}

\adbb{
Returning to~\Ref{ADB-Tiffany-2}, we note that
\eqs
    \lefteqn{ \lti \pr_{x_0}[\{t^{-1/2}(\log X_t - \bb t) \ge y + t^{-1/2}r(t)\} \cap \cX] }\\
     &&\Eq \lti \pr_{x_0}[\{t^{-1/2}(\log X_t - \bb t) \ge y - t^{-1/2}r(t)\} \cap \cX]
     \Eq \pr_{x_0}[\cX](1 - \Phi(y/v)),
\ens
by \Ref{ADB-last-line} and because $t^{-1/2}r(t) \to 0$ as $t \to \infty$.  It thus follows that
\[
     \limsup_{t\to\infty}|\pr_{x_0}[\{t^{-1/2}(\log X_{t+\D(t)} - \bb t) \ge y\} \cap \cX] - \pr_{x_0}[\cX](1 - \Phi(y/v))|
         \Le \limsup_{t\to\infty} \pr_{x_0}[\cE_{s,t}^c \cap \cX]
\]
for any $s \ge 1$, and hence, from \Ref{ADB-Xc-prob} and~\Ref{ADB-limlim}, that
\[
    \lti \pr_{x_0}[t^{-1/2}(\log \Yh_t - \bb t) \ge y] \Eq \pr_{x_0}[\cX](1 - \Phi(y/v))
\]
for all $y \in \re$.  The conclusion of the theorem now follows by replacing~$t$ by $h_0(\snt) - h_0(\snt_0)$ for 
$\snt \ge \snt_0$,
$\snt \in \nat$, and noting that $h_0(\snt) - h_0(\snt_0) \sim \log(\snt/\snt_0)$ as $\snt \to \infty$. \ep
}
\end{proof}

\subsection{The discrete degree weighted process}
\adbb{
There is also a degree weighted process~$\tY$ in discrete time.  Multiplying equation~\Ref{ADB-probability-equation}
for~$\pp_{\snt+1,k}$ by~$k$, it follows that, with $u_{\snt,k} := k\pp_{\snt,k}$,
\eqs
    u_{\snt+1,k} &=&   u_{\snt,k} + \frac1{\snt+1}\,\Bigl\{ - u_{\snt,k}(1 + \a k) + \a k u_{\snt,k-1}  
    + q_k u_{\snt,k}  +  \sum_{j \ge k} p_j(1 - q_j)u_{\snt,j} \tpi_{jk} \Bigr\},  
\ens
where $\tpi_{jk} := k\p_{jk}/\{jp_j\}$, so that~$\tPi_j$ is the size--biased transformation of~$\Pi_j$.
This can be re-written in vector form as
\[
    u_{\snt+1}\TT \Eq u_\snt\TT\Bigl\{I + \frac1{\snt+1}\,\{[\tQ]_{\snt+1} + (2\a-1)I\} \Bigr\},
\]
where~$\tQ$ is the $Q$-matrix defined in~\Ref{wtQmat} with $\b_j := 1-q_j$; recall that, for each $j \ge 1$, 
$0 < q_j + p_j(1-q_j) = \a < 1$
is constant.  Hence, defining $v_{\snt} := \prod_{s = \snt_0}^{\snt-1}\{1 + (2\a-1)/(s+1)\}^{-1} u_{\snt}$, this implies that,
in parallel to~\Ref{ADB-Q-version-discrete},
\[
    v_{\snt+1}\TT \Eq v_\snt\TT\Bigl\{I + \frac1{\snt+2\a}[\tQ]_{\snt+1}\Bigr\}, \quad \mbox{and so}\quad
           v_\snt\TT \Eq  v_{\snt_0}\TT \prod_{j=\snt_0+1}^{\snt} \Bigl\{I + \frac1{j+2\a-1}[\tQ]_j\Bigr\}.
\]
These equations show that, if the initial configuration at time~$\snt_0$ has~$j_0$ vertices, then the family of 
vectors~$j_0^{-1}(v_\snt,\,\snt\ge \snt_0)$ can be interpreted as the sequence of probability 
distributions~$(\tpp_\snt\ujoto,\,\snt\ge \snt_0)$ of a time inhomogeneous
Markov chain~$\tY$ starting in state~$j_0$ at time~$\snt_0$, since $v_{\snt_0} = u_{\snt_0} = j_0$ if $\pp_{\snt_0}$ 
is the distribution
with point mass at~$j_0$, and since $I + (j+2\a-1)^{-1}[\tQ]_{j}$ is a stochastic matrix, for each~$j$.
This in turn implies that, if $\pp_{\snt}\ujoto$ denotes the distribution of~$Y$, starting in state~$j_0$ at time~$\snt_0$, 
then 
\[
       \pp_{\snt,k}\ujoto \Eq j_0 k^{-1} \prod_{s = \snt_0}^{\snt-1}\{1 + (2\a-1)/(s+1)\} \tpp_{\snt,k}\ujoto,
\]
so that the fractions~$\pp_{\snt}\ujoto$ can be deduced from the probabilities~$\tpp_t\ujoto$.  This is the discrete analogue
of Proposition~\ref{quasi}; note that, for large~$\snt_0$ and $\snt \ge \snt_0$, 
\[
         \prod_{s = \snt_0}^{\snt-1}\{1 + (2\a-1)/(s+1)\} \ \approx\ e^{-(1-2\a)(\snt-\snt_0)},
\]
making the comparison clearer.
} 

\adbb{
A time changed version~$\tYh$ of
the discrete time inhomogeneous degree weighted process~$\tY$ can be closely coupled to the continuous
time degree weighted process, using arguments exactly like those above.  The processes $\tY$ and~$\tX$ have the same 
jump chains,
and can be constructed in parallel, using the same realizations of the jump chain, and the same sequence of independent 
uniform random variables.  The time scale for the discrete process that matches that of the continuous time
process is $(h_a(\snt) - h_a(\snt_0),\, \snt \ge \snt_0)$, where $a := 2\a-1$, and the waiting time until 
the next jump in the time scaled discrete process~$\tYh$, when starting in state~$k$ at time $h_a(\snt) - h_a(\snt_0)$
(or at step~$\snt$ of the process~$\tY$), has the distribution of~$V_a(\snt,-\tQ_{k,k})$.
The subsequent arguments, showing that, with this coupling, $\tYh_t = \tX_{t + \tDe(t)}$
for all $t \ge 0$, for a difference~$\tDe(t)$ defined
analogously to the definition~\Ref{ADB-Delta-def} of~$\D(t)$, and that $\tDe(t)$ converges a.s.\ as $t\to\infty$,
are substantially the same.
}

\adbb{
The next lemma can be applied to show that the limiting probabilities for the degree weighted inhomogeneous chain~$\tY$
are the same as those for~$\tX$,
if the process~$\tX$ is positive recurrent.  The lemma can also be applied to the discrete versions of
the unweighted chain, in models, such as the DD model with `random re-wiring' in Section~\ref{ADB-sect-rewire}, 
in which the original chain is positive recurrent.  For the basic DD model, with~$0$ an absorbing state, it is not needed.
}

\begin{lemma}\label{ADB-convergence-equivalence}
\adbb{
 Suppose that $(X_t,\,t\ge0)$ is an irreducible positive recurrent pure jump Markov process on~$\Z$, 
with some initial distribution~$\l$ and with stationary distribution~$\pi$.  Let
$(\D(t),\,t\ge0)$ be a c\`adl\`ag stochastic process such that $\D := \lti \D(t)$ exists
and is finite a.s.  Then, for each $j \in \Z$,
\[
     \lti \pr[X_{t+\D(t)} = j] \Eq \pi_j.
\]
}
\end{lemma}

\begin{proof}
\adbb{
 Suppose first that $\D$ is integrable.  Write $\hDe(t) := \ex\{\D \giv \cF_t\}$, where
$\cF_t := \sigma\{X_s,\,0\le s\le t\}$.  Define the events
\[
     A_s(\eps) \Def \{|\D(t) - \D| \le \eps \mbox{ for all }t \ge s\};\quad
     \hA_s(\eps) \Def \{|\hDe(t) - \D| \le \eps \mbox{ for all }t \ge s\},
\]
noting that, for each $\eps > 0$,
\[
      \lim_{s\to\infty} \pr[A_s(\eps)] \Eq \lim_{s\to\infty} \pr[\hA_s(\eps)] \Eq 1.
\]
}

\adbb{
Now, for any $j \in \Z$, $s,t > 0$ and $\eps > 0$, 
\eqs
   \{X_{t+\D(t)} = j\} &\supset& \{X_{t + \hDe(t) + u} = j,\,-\eps \le u \le \eps\}\cap A_s(\eps/2) \cap \hA_s(\eps/2) \\
             &\supset& \{X_{t + \hDe(s) + u} = j,\,-\eps \le u \le \eps\}\cap A_s(\eps/2) \cap \hA_s(\eps/4),
\ens
from which it follows that
\eq\label{ADB-prob-lb}
    \pr[X_{t+\D(t)} = j] \ \ge\ \pr[X_{t + \hDe(s) + u} = j,\,-\eps \le u \le \eps]  
                       - \pr[A_s^c(\eps/2)] - \pr[\hA_s^c(\eps/4)].
\en
Then, by the Markov property,
\eqs
     \lefteqn{ \pr[X_{t + \hDe(s) + u} = j,\,-\eps \le u \le \eps \giv \cF_s]\,I[t + \hDe(s) - \eps \ge s] }\\
         &&\Eq e^{-2q_j\eps}\, \pr[X_{t + \hDe(s) - \e} = j \giv \cF_s]\,I[t + \hDe(s) - \eps \ge s] \\
         &&\Eq e^{-2q_j\eps}\,P_{X_s,j}(t + \hDe(s) - \e - s)\,I[t + \hDe(s) - \eps \ge s],
\ens
where~$q_j$ denotes the jump rate from the state~$j$, and $P_{i,j}(u) := \pr[X_u=j \giv X_0 = i]$.
Hence, and because $\lti I[t + \hDe(s) - \eps \ge s] = 1$ a.s., we have
\eqs
   \lefteqn{ \lti \pr[X_{t + \hDe(s) + u} = j,\,-\eps \le u \le \eps \giv \cF_s] }\\
  &&\Eq
    \lti\{\pr[X_{t + \hDe(s) + u} = j,\,-\eps \le u \le \eps \giv \cF_s]\,I[t + \hDe(s) - \eps \ge s]\}
                  \\ 
   &&\Eq e^{-2q_j\eps}\, \pi_j\quad \mbox{a.s.}\,,
\ens
and thus, by dominated convergence,
\[
    \lti\pr[X_{t + \hDe(s) + u} = j,\,-\eps \le u \le \eps] \Eq e^{-2q_j\eps}\, \pi_j.
\]
It now follows from~\Ref{ADB-prob-lb} that
\[
    \liminf_{t \to \infty} \pr[X_{t+\D(t)} = j] \ \ge\ e^{-2q_j\eps}\, \pi_j - \pr[A_s^c(\eps/2)] - \pr[\hA_s^c(\eps/4)],
\]
for any $s > 0$, $j \in \Z$ and $\eps > 0$.  Letting $s \to \infty$ and then $\eps \to 0$ thus implies that
\eq\label{ADB-individual-lb}
    \liminf_{t \to \infty} \pr[X_{t+\D(t)} = j] \ \ge\ \pi_j
\en
for all $j \in \Z$.  The inequality~\Ref{ADB-individual-lb}, used in the second and final inequalities below, 
now implies that, for any $j_0 \in \Z$  and any {\it finite\/} set $J \subset \Z$ such that $j_0 \in \Z$,
\eqa
     1 - \sum_{j \in J} \pi_j &\ge& 
   \limsup_{t \to \infty} \pr[X_{t+\D(t)} \in J] - \sum_{j \in J} \pi_j \label{ADB-first-line}\\
     &\ge& 
       \limsup_{t \to \infty} \,(\pr[X_{t+\D(t)} = j_0] - \pi_{j_0}) 
    \ \ge\ 
           \liminf_{t \to \infty} \,(\pr[X_{t+\D(t)} = j_0] - \pi_{j_0}) \ \ge\ 0. \non
\ena
Letting the set~$J$ increase towards the whole of~$\Z$, the sum $\sum_{j \in J} \pi_j$ can be made arbitrarily
close to~$1$, and so the left hand side of~\Ref{ADB-first-line} can be made arbitrarily small, implying that
\[
    \lim_{t \to \infty} \pr[X_{t+\D(t)} = j_0] \Eq \pi_{j_0}
\]
for any $j_0 \in \Z$, as required.
}

\adbb{
If $\D$ is not integrable, for any $M > 0$, define $\D_M := [\D]_{-M}^M := \max\{-M,\min(\D,M)\}$ to be $\D$ truncated
within the range $[-M,M]$, and let $\D_M(t) := [\D(t)]_{-M}^M$.  Then, because $[\cdot]_{-M}^M$ is
a continuous function, $\lti \D_M(t) = \D_M$ a.s., and hence, by the previous result,
\[
     \lim_{t \to \infty} \pr[X_{t+\D_M(t)} = j] \Eq \pi_{j},\qquad j \in \Z,
\]
since $\D_M$ is integrable.   But
\[
   |\pr[X_{t+\D_M(t)} = j] - \pr[X_{t+\D(t)} = j]| \Le \pr[\D_M(t) \neq \D(t)] \Eq \pr[|\D(t)| > M],
\]
and, since $\D(t) \to \D$ a.s., $\limsup_{t \to \infty} \pr[|\D(t)| \ge M] \Le \pr[|\D| \ge M]$.
Hence, for all $M > 0$,
\[
   \p_j - \pr[|\D| \ge M] \Le \liminf_{t\to\infty} \pr[X_{t+\D(t)} = j] \Le
                 \limsup_{t\to\infty} \pr[X_{t+\D(t)} = j] \Le \p_j + \pr[|\D| \ge M],
\]
implying that $\lti\pr[X_{t+\D(t)} = j] = \p_j$, and completing the proof. \ep
}
\end{proof}

\adbb{
Applying Lemma~\ref{ADB-convergence-equivalence} to the process~$\tYh$ shows that, if~$\tX$ is positive recurrent,
with stationary distribution~$\tp(\infty)$, then, in view of~\Ref{ADB-matching-processes},
\[
      \lti \pr[\tYh_t = j] \Eq \tp_j(\infty), \qquad j \ge 1.
\]
It thus also follows that 
\[
     \lim_{\snt \to \infty} \pr[\tY_{\snt} = j] \Eq \tp_j(\infty), \qquad j \ge 1,
\]
and hence that the inhomogeneous process~$\tY$ has the same limiting distribution as~$\tX$,
since $\tY_{\snt} = \tYh_{h_0(\snt) - h_0(\snt_0)}$, and $\lim_{\snt\to\infty} \{h_0(\snt) - h_0(\snt_0)\} = \infty$.
}

\section{Variants of the basic process}\label{ADB-immigration}
\setcounter{equation}{0}

\subsection{The process with removal of edges}\label{ADB-sect-deaths}
 A process, in which edges are also allowed to be removed, has been discussed in \cite{hermann2019partial}.
The analogue in
 our formulation is obtained by modifying the process~$X$ of Section~\ref{ADB-absorption} by allowing deaths at 
rate $k\d_k$, $k\ge1$, with $|\d_k - \d| \le c_5 k^{-\g}$ for all~$k$.  This gives a $Q$-matrix
\begin{align}
   &Q_{k,k+1}\ud \Eq k\a_k ;\qquad Q_{k,k-1}\ud \Eq k\d_k;\qquad Q_{k,k}\ud \Eq -\{k(\a_k + \d_k) + \bet_k(1-\p_{kk})\}; \non\\
   &Q_{k,j}\ud \Eq  \bet_k\p_{kj},\quad 0\le j\le k-1; \qquad Q_{k,j}\ud \Eq 0,\quad j \ge k+2,  \label{ADB-Q-matrix-d}
\end{align}
The properties of the modified process~$X\ud$ are much as for the process without deaths.  First, the asymptotic 
growth rate for the process without catastrophes becomes $\a\ud := \a - \d$, rather that~$\a$.
This affects the criterion for transience and recurrence in the obvious way;  Theorem~\ref{ADB-recurrence} holds, 
but with $\a$
replaced by~$\a\ud$, and the proofs are otherwise unchanged.  The net growth rate for the process with catastrophes 
now becomes $\bb\ud := \a\ud - \bet\log(1/p)$,
and Theorem~\ref{aslimit} holds with $\a\ud$ for~$\a$ and with $\bb\ud$ for~$\bb$.  Here, some small modifications to the
proofs of two of the lemmas are needed, largely because the process is no longer increasing between catastrophes.

    First, in the proof of Lemma~\ref{ADB-A-and-B-probs}, the martingale~$M$ in~\Ref{ADB-M-def-2} has to be modified.
Here, $\a_{Y\uii_u}$ is to be replaced by~$\tal_{Y\uii_u}$, where, in order to obtain a martingale, we define 
$\tal_j := \a_j - j\d_j/(j-1)$ in $j \ge 2$; the quantity~$\tal_j$ is not quite the same as~$\a_j - \d_j$, but
is close to it when~$j$ is large, which is the case in the arguments of Lemma~\ref{ADB-A-and-B-probs}.  A stopped
version~$M^\s$ is used in the proofs, where $M^\s_l := M_{l\wedge\s}$, and 
\[
    \s \Def \s\uii \Def \inf\bigl\{ l \ge 1\colon Y\uii_{\Xss{i-1}+l} \le \half\{\Xss{i-1} + l\r(\half\Xss{i-1})\} \bigr\},
\]
where $\r(j) := \inf_{k \ge j}\{(\a_k - \d_k)/(\a_k + \d_k)\}$, a lower bound for the drift
of the asymmetric random walk described by the jump chain of~$Y\uii$ in $\Z\cap [j,\infty)$. Provided that
$\Xss{i-1} \ge 2$, there thus is no need to define~$\tal_j$ for $j \le 1$.  It is then possible to
apply Doob's inequality to control the fluctuations of~$M^\s$, and this can be translated into bounds
on the probability of the events~$A_i$, again with $\a$ replaced by~$\a\ud$.
Note that the variance of~$M^\s$ and the means and
variances of $\D^\s_l := \D_{l\wedge\s}$, conditional on~$\cF_{T_{i-1}}$, can be bounded much as in \Ref{ADB-M-variance} and
\Ref{ADB-Delta-mean-var-bnd}, but with different constants.  For instance, temporarily defining
\[
      \hp_j \Def 1 - \hq_j \Def \frac{\a_j}{\a_j + \d_j},
\]
we have
\eqa
    \lefteqn{ \ex\bigl\{(M_l - M_{l-1})^2 \giv \cF_{S_{l-1}} \cap \{\s\uii > l-1\} \cap \{Y\uii_{S_{l-1}} = j\} \bigr\} }\non\\
            &&\Eq \hp_j\hq_j\Bigl\{\frac1j + \frac1{j-1}\Bigr\}^2 + \Bigl\{\frac{\hp_j}j - \frac{\hq_j}{j-1}\Bigr\}^2 
            \Le 4j^{-2},  \label{ADB-MG-var-new}
\ena
in $j \ge 2$, where the first term in~\Ref{ADB-MG-var-new} comes from the variability in $h(Y\uii_{S_l})$ and the
second from the variability in $\tal(S_l - S_{l-1})$.  Hence, for any $l \ge 1$
\eq\label{ADB-variance-sum-new}
     \var\{M^\s_l \giv \cF_{T_{i-1}}\} \Le 4 \sum_{k \ge 1} \ex\{(Y\uii_{S_{k-1}})^{-2} I[\s\uii > k-1] \giv \cF_{T_{i-1}}\}.
\en
But since, for any $r > 1$,
\eqs
      \sum_{l \ge 0} (Y\uii_{S_l})^{-r} I[\s\uii > l] 
           &\le& \sum_{l \ge 0} \bigl\{ \half\{\Xss{i-1} + l\r(\half\Xss{i-1})\} \bigr\}^{-r} \\
           &\le&  \frac{2^r}{\Xss{i-1}^{r-1}}\,\Bigl\{\frac1{\Xss{i-1}} + \frac{1}{(r-1)\r(\half\Xss{i-1})}\Bigr\},
\ens
the sum in~\Ref{ADB-variance-sum-new} is bounded by $16 \Xss{i-1}^{-1} (1 + 1/\r(j_0))$,
uniformly in $\Xss{i-1} \ge 2j_0$, for any~$j_0$ such that $\r(j_0) > 0$. This is the analogue of~\Ref{ADB-M-variance}.
The bounds analogous to \Ref{ADB-M-variance} and \Ref{ADB-Delta-mean-var-bnd}
in turn imply analogues of the probability bounds \Ref{ADB-M-bound} and~\Ref{ADB-Delta-bound}
for the deviations of $M^\s$ and~$\D^\s$, again with different constants.
By using a Wald martingale for the simple
random walk with drift~$\r(\half\Xss{i-1})$, the probability $\pr[\s\uii < \infty \giv \cF_{T_{i-1}}]$ is 
bounded by $e^{-c\Xss{i-1}}$,
for a suitable constant~$c = c(\half\r(\Xss{i-1}))$.   This is enough to establish Lemma~\ref{ADB-A-and-B-probs}, under
the extended model.

    In the proof of Lemma~\ref{ADB-first-growth}, the fact that the process~$Y\uii$ is non-decreasing is used, 
but only for convenience.
The proof can easily be modified for the process in which deaths are allowed to occur.
It is enough, instead, to know that~$Y\uii$ does not fall below the value $\half\Xss{i-1}$, and to  
define~$\h$ so that $(2/\h)^{4/\g} > 2j'$ in the analogue of Lemma~\ref{ADB-first-growth}.  Then, in the
induction argument used to prove Lemma~\ref{ADB-first-growth}, note that on~$A_i$, for any $t \in [T_{i-1},T_i)$,
\[
    \log X_t \ \ge\ \log \Xss{i-1} - \Xss{i-1}^{-\g/4} \ \ge\ \log(\Xss{i-1}/2),
\]
as desired, if $\h < 2\log 2$, by the induction hypothesis, so that then $\bet_{X_t} \le \bet'$ for all
$T_{i-1} \le t < T_i$.  The remainder of the proof is as before.
In the analogue of Lemma~\ref{ADB-prob-HE}, $\h$ should also now satisfy $(2/\h)^{4/\g} > 2j'$.

Theorem~\ref{CLT} also holds, with $\bb\ud$ for~$\bb$.  As for Theorem~\ref{aslimit}, some small modifications in the
proofs are needed, again because the process is no longer increasing between catastrophes.

The biological motivation for allowing edges to be removed is that some interactions may decline in importance
over time, being replaced by more advantageous interactions.  However, these removals would not occur only at
the occasions on which a vertex is duplicated, so that the original discrete process is not well suited to
such a modification.  \cite{hermann2019partial} introduce edge deletions in a process very similar to the
continuous time model of \cite{jordan}, in which edges can easily be modelled as having independent exponential
lifetimes.

In \cite{hermann2020markov}, extra terms are introduced into the birth--death--catastrophe process, allowing for upward
jumps of sizes greater than one, those corresponding to a Markov branching process.  Here, we can again allow such jumps,
now allowing the jump rates to vary a little with the state~$k$. Writing $\nat_{-1} := \{-1\}\cup\nat$, we can modify 
the matrix~$Q$ further by setting
\begin{align}
   &Q_{k,k+j}\ub \Eq k a_{k,j},\quad j \in \nat;\qquad Q_{k,k-1}\ub \Eq k a_{k,-1} + \bet_k\p_{k,k-1};\non\\
   &Q_{k,k}\ub \Eq -\Bigl\{k\sum_{j\in\nat_{-1}} a_{k,j} + \bet_k(1-\p_{kk}) \Bigr\}; \qquad
   Q_{k,j}\ub \Eq  \bet_k\p_{kj},\quad 0\le j\le k-2.  \label{ADB-Q-matrix-b}
\end{align}
Here, we assume that $\inf_{k \ge 1} \sum_{j\in\nat_{-1}} a_{k,j} > 0$, and that
\[
    |a_{k,j} - a_{*,j}| \Le  c_{*,j} k^{-\g}, \quad j \in \nat_{-1},\quad k \ge 1,
\]
where 
\[
    \sum_{j \in \nat_{-1}} j a_{*,j}\ =:\ \a\ub \ <\ \infty;\qquad \sum_{j \in \nat_{-1}} jc_{*,j} < \infty,
\]
and $\sum_{j \in \nat_{-1}} j^2 a_{k,j}$ is uniformly bounded in~$k$.  It then follows that the asymptotic 
growth rate for the process without catastrophes is $\a\ub$, rather that~$\a$, which is consistent with~$\a\ud$
when $a_{*,j} = 0$ for all $j \ge 2$, and with~$\a$ in Section~\ref{ADB-absorption}, when also $a_{*,-1} = 0$.
Largely using the proof of Theorem~\ref{ADB-hX-recurrence}, we can show that the process~$\sXb$, that has an
additional transition from state~$0$ to state~$1$, satisfies the following theorem.

\begin{theorem}\label{ADB-sXb-recurrence}
For the process~$\sXb$, we have the following behaviour:
\begin{align}
 &\mathrm{(1)}:\quad \mbox{if}\ \alpha\ub < \bet \log(1/p),\ \mbox{then~$\sXb$ is geometrically ergodic};\non\\
 &\mathrm{(2)}:\quad \mbox{if}\ \alpha\ub = \bet \log(1/p),\ \mbox{then~$\sXb$ is null recurrent};\label{ADB-FLT-hat-X-d}\\
 &\mathrm{(3)}:\quad \mbox{if}\ \alpha\ub > \bet \log(1/p),\ \mbox{then~$\sXb$ is transient}.\non
\end{align}
\end{theorem}

\begin{proof}
 The difference from the proof of Theorem~\ref{ADB-hX-recurrence} comes from needing to find the leading terms in
$k \sum_{j \in \nat_{-1}} a_{k,j}(f(k+j) - f(k))$, rather than just in $k(f(k+1) - f(k))$, for the same four choices
of~$f$.  For the two recurrence arguments, when~$f$ is concave, we immediately have
\[
    \sum_{j \in \nat} a_{k,j} (f(k+j) - f(k)) \Le \sum_{j \in \nat} ja_{k,j} f'(k),
\]
and the rest of the proof is straightforward, under the assumptions on the~$a_{k,j}$.  For the remaining two arguments,
we use~\Ref{ADB-Jensen-lower} with $x = k$ and $y = k+j$ in proving null recurrence, and~\Ref{ADB-Jensen-upper}
in proving transience. \ep
\end{proof}

\nin The analogue of Theorem~\ref{ADB-recurrence-2} is also true, in case~(1), with~$\a\ub$ in place of~$\a$.

In case~(3), the analogue of Theorem~\ref{aslimit}, with~$\a\ub$ in place of~$\a$, is also true.  As discussed above
for~$X\ud$, there are modifications to the proofs of Lemmas \ref{ADB-A-and-B-probs}, \ref{ADB-first-growth}
and~\ref{ADB-prob-HE}, because of the deaths.  In Lemma~\ref{ADB-A-and-B-probs}, the martingale
$M$ in~\Ref{ADB-M-def-2} again has to be modified, now with $\a_{Y\uii_u}$ replaced by~$\tal_{Y\uii_u}\ub$, where we define 
\[
    \tal_j\ub \Def j\sum_{k \in \nat_{-1}}  a_{j,k} (h(j+k-1) - h(j-1))\ \approx\ \a\ub
\]
in $j \ge 2$.  The lower bound~$\r(j)$ on the drift of the jump chain, used in defining the stopping time~$\s\uii$,
also needs to be modified to
\[
    \r(j) \Def \inf_{k \ge j}\Bigl\{\frac{\sum_{l \in \nat_{-1}} l a_{k,l}}{\sum_{l \in \nat_{-1}} a_{k,l}} \Bigr\}\,.
\]
As for~$X\ud$, it is necessary to be able to bound the variance of~$M^\s$ and the means and
variances of $\D^\s_l := \D_{l\wedge\s}$, conditional on~$\cF_{T_{i-1}}$, much as in \Ref{ADB-M-variance} and
\Ref{ADB-Delta-mean-var-bnd}, but with different constants.  Under our assumptions on the values of~$a_{k,j}$,
the calculations are routine.  Finally, again with small modifications to the proof, Theorem~\ref{CLT} also holds 
for~$X\ub$, with $\bb\ub$ for~$\bb$.

\subsection{The process with random re-wiring}\label{ADB-sect-rewire}

In \cite{pastor2003evolving}, 
the basic DD model is modified to accommodate random connections,
added to the duplicate vertex after copying and thinning.  At time~$(\snt+1)$, each vertex, other than the vertex being
copied and its neighbours, may be connected to the duplicate, independently with probability~$r/\snt$.

In this case, the equations~\Ref{ADB-probability-equation} for probabilities in the basic DD model, for $k \ge 1$, 
are modified to
\eqa
    \pp_{\snt+1,k}\ur &=&   \pp_{\snt,k}\ur + \frac1{\snt+1}\,\Bigl\{ - \{1 + \a k + r(1-(k+1)/\snt)\}\pp_{\snt,k}\ur  
           \non\\
      &&\qquad\mbox{} + \{(k-1)\a + r(1 - k/\snt)\} \pp_{\snt,k-1}\ur 
               \label{ADB-probability-equation-rewire} \\
       &&\qquad\qquad\mbox{}   + q \sum_{j \le k}\pp_{\snt,j}\ur \hpi_{jk}\urt
          +  \sum_{j \ge 0} (1 - q)\pp_{\snt,j}\ur \p_{jk}\urt \Bigr\}, \non
\ena
where, as before, $\a := q + p(1-q)$; the probability distribution~$\hPi_k\urt$ is the binomial 
distribution~$\Bi(\snt-k,r/\snt)$, having point probabilities~$\hpi_{jk}\urt$, and the probability 
distribution~$\Pi_k\urt$, with point probabilities~$\p_{jk}\urt$, is the convolution of $\Pi_k$ and~$\hPi_k\urt$, 
and has mean~$kp + r - rk/\snt$.  In vector notation, the equations can be written as
\eq\label{ADB-Q-version-discrete-rewire-0}
     (\pp_{\snt+1}\ur)\TT \Eq (\pp_{\snt}\ur)\TT\{I + (\snt+1)^{-1}[Q^{(r,\snt+1)}]_{\snt+1}\}, 
\en
and so
\eq\label{ADB-Q-version-discrete-rewire}
          (\pp_\snt\ur)\TT \Eq (\pp_{\snt_0}\ur)\TT \prod_{j=\snt_0+1}^\snt \{I + j^{-1}[Q^{(r,j)}]_j\},
\en
where
\begin{align}
   &Q_{k,k+1}\urt \Eq \a k + r(1 - (k+1)/\snt) + (1-q)\p_{k.k+1}\urt + q\hpi_{k,k+1}\urt;\non\\
   &Q_{k,k}\urt \Eq -\{\a k + r(1 - (k+1)/\snt) + (1-q)(1-\p_{kk}\urt) + q(1 - \hpi_{kk})\urt\}; \label{ADB-Q-matrix-rewire}\\
   &Q_{k,j}\urt \Eq  (1-q)\p_{kj}\urt,\quad 0\le j\le k-1;\quad 
          Q_{k,j}\urt \Eq (1-q)\p_{kj}\urt + q\hpi_{k,j}\urt, \quad j \ge k+2,\non              
\end{align}
for each $k \ge 1$, and with 
\eqa
    &&Q_{0,1}\urt \Eq r(1 - 1/\snt) + \hpi_{01}\urt;\qquad Q_{0,j}\urt \Eq \hpi_{01}\urt, \quad j \ge 2; \non\\
    &&Q_{0,0}\urt \Eq -\{r(1 - 1/\snt) + 1 - \hpi_{00}\urt\}.  \label{ADB-Q-matrix-rewire-0}
\ena
Note that, if $r=0$, $Q\urt$ reduces to the $Q$-matrix for the basic DD model. If~$\snt$ is large, the $Q$-matrices 
$Q\urt$ are very close to the matrix~$Q\ur$, in which, for all $k \ge 0$, the fractions $(k+1)/\snt$ are replaced by zero,
$\hPi_k\urt$ is replaced by~$\hPi\ur$, the Poisson distribution with mean~$r$, and~$\Pi_k\urt$ is replaced by~$\Pi_k\ur$,
the convolution of~$\Pi_k$ and~$\hPi\ur$.
 
Now the process~$X\ur$ with $Q$-matrix $Q\ur$ is almost of the same form as that for the irreducible process~$\sX$ in
Section~\ref{ADB-absorption}, with $\a_0 := r$, and with $\a_k := \a + r/k$ and $\bet_k := (1-q)$ for 
all $k \ge 1$. The difference is that it also has extra elements in the rates, because of the (limited) re-wiring.

First, instead of having rate~$\a k$ of upward jumps in state~$k$, with~$\a$ constant, we have
a rate~$k a_{k,1}\ur$, where
\eq\label{ADB-immigration-alpha}
   a_{k,1}\ur \Eq \a + \{r + (1-q)\p_{k.k+1}\ur + q\hpi_{k,k+1}\ur\}/k \Eq \a + O(k^{-1}), \quad k\ge1,
\en
and there are also jump rates $k a_{k,j}\ur$ for $j \ge 2$ as in~\Ref{ADB-Q-matrix-b}, with
\[
     a_{k,j}\ur \Eq \{(1-q)\p_{k.k+j}\ur + q\hpi_{k,k+j}\ur\}/k \Eq O(k^{-1}), \quad k \ge 1.
\]
Thus the framework is technically that of the process~$X\ub$, as given by~\Ref{ADB-Q-matrix-b},
though the quantities~$a_{*,j}$ are all zero, except for $j = 1$.
Then, instead of the copy distribution from a vertex of degree~$k$ being a mixture of the point mass on~$k$
with probability~$q$ and the distribution~$\Pi_k$ with probability~$1-q$, the distribution is convolved with
the Poisson distribution~$\Po(r)$.  This in part contributes to the extra terms~$a_{k,j}\ur$, but also has the consequence
that the distribution~$\Pi_k\ur$, when restricted to the set $\{0,1,\ldots,k\}$, has probability mass less than~$1$.
In order to fit with a process of the form~$X\ub$, the factor~$1-q$ should be replaced by a parameter~$\b_k\ur$,
where $\b_k\ur := (1-q)\Pi_k\ur\{[0,k]\}$, and the downward jump distribution should be modified to 
$\tPi_k\ur := (1-q)\Pi_k\ur/\b_k\ur$.
Under the assumptions on~$\Pi_k$, $\b_k\ur = (1-q) + O(k^{-\g})$, the mean of~$\tPi_k\ur$ is asymptotically 
$kp(1 + O(k^{-\g}))$, and its variance is of order~$O(k^{2-\g})$.

Thus the process~$X\ur$ can be analyzed in the same way as for~$X\ub$.  
For instance, Theorem~\ref{ADB-sXb-recurrence} shows that~$X\ur$ is positive recurrent if  $\alpha < (1-q) \log(1/p)$,
null recurrent if   $\alpha = (1-q) \log(1/p)$ and transient if   $\alpha > (1-q) \log(1/p)$.
In the case of positive recurrence, Theorem~\ref{ADB-recurrence-2} indicates which moments of the stationary
distribution exist.  Note that the critical exponent~$\h^*$ is smaller than~$(1-q)/\a = (1-q)/(q + p(1-q)) \le 1/p$,
so that higher moments can only be finite when~$p$ is small enough.  In the case of transience, as for the
process~$X\ub$, the conclusions of Theorems \ref{aslimit} and~\ref{CLT} hold.

We now show that the behaviour of the inhomogeneous Markov chain~$\hY\ur$, whose probabilities are governed
by the equations~\Ref{ADB-Q-version-discrete-rewire}, can indeed be deduced from that of the process~$X\ur$.
To do this, we couple~$\hY\ur$ to a process~$Y\ur$ whose probabilities are given by~\Ref{ADB-Q-version-discrete-rewire},
but with~$Q^{(r,j)}$ replaced by~$Q\ur$.  If the two processes are started together in the state $(\snt_1,j_1)$, 
construct them jointly so as to make identical transitions as far as possible.  If, at any step~$\snt$, the two processes 
are in the same state~$k$, the jump probabilities differ in two respects.  First,
the probability of jumping up by~$1$ because of a re-wiring from a copied vertex that is not a neighbour
is an amount $r(k+1)/\{\snt(\snt+1)\}$ larger in $Y\ur$ than in~$\hY\ur$.  Secondly, the copy distributions differ,
in that one of them is~$\Pi_k$ convolved with $\Bi(\snt - (k+1),r/\snt)$ and the other with~$\Po(r)$.
As in \cite{prokhorov1953asymptotic}, the total variation distance between these two distributions is at most
$c\snt^{-1}$, where~$c$ can be taken to be~$r$.  Hence the jump distributions of $Y\ur$ and~$\hY\ur$ when leaving state~$k$
at time~$\snt+1$ differ in total variation by at most $r(k+2)/\{\snt(\snt+1)\}$.
This suggests that the two processes can be coupled so as to have identical paths with high
probability, if~$\snt_1$ is large enough.  

More precisely, the probability that the two processes fail to remain identically coupled is at most
\[
    E_{j_1,\snt_1} \Eq  r\sum_{\snt > \snt_1} \frac{\ex\{Y\ur_{\snt}+2 \giv Y\ur_{\snt_1} = j_1\}}{\snt(\snt+1)}.
\]
Now, from the one-step probabilities for~$Y\ur$, we have
\[
     \ex\{Y\ur_{\snt+1} \giv Y\ur_{\snt} = j\} \Le j + \frac{\a j + 2r + (1-q)(p-1)j}{\snt+1} 
                  \Eq j\Bigl\{1+\frac{2\a-1}{\snt+1}\Bigr\} + \frac {2r}{\snt+1},
\]
and hence, for $\snt \ge \snt_1$, 
\[
   \ex\{Y\ur_{\snt+1} \giv Y\ur_{\snt_1} = j_1\} 
       \Eq \Bigl\{1+\frac{2\a-1}{\snt+1}\Bigr\}\ex\{Y\ur_{\snt} \giv Y\ur_{\snt_1} = j_1\} + \frac {2r}{\snt+1}\,.
\]
Iterating, we deduce that
\[
    \ex\{Y\ur_{\snt} \giv Y\ur_{\snt_1} = j_1\}  \Le 
     \begin{cases}
            \Bigl(j_1 + \frac{2r}{2\a-1}\Bigr) \prod_{l = \snt_1+1}^{\snt} \Bigl\{1 + \frac {2\a-1} l\Bigr\},
                                  &\mbox{ if } 2\a > 1;\\
            j_1 +  \sum_{l = \snt_1+1}^{\snt} \frac{2r}l\,,  &\mbox{ if } 2\a = 1;\\
            j_1\prod_{l = \snt_1+1}^{\snt} \Bigl\{1 + \frac {2\a-1} l\Bigr\} + \frac{2r}{1-2\a}\,, &\mbox{ if } 2\a < 1 ,
     \end{cases}
\]
 for $\snt > \snt_1$, and hence that 
\eq\label{ADB-Yur-mean}
    \ex\{Y\ur_{\snt} \giv Y\ur_{\snt_1} = j_1\} \Le 
     \begin{cases}
            C(j_1 + 2r/(2\a-1))(\snt / \snt_1)^{2\a-1},   &\mbox{ if } 2\a > 1;\\
            j_1 + 2r \{\log(\snt / \snt_1) + \g\},   &\mbox{ if } 2\a = 1;\\    
            C j_1 ( \snt_1/\snt)^{1-2\a} + 2r/(1-2\a),  &\mbox{ if } 2\a < 1:
     \end{cases}
\en
here, $\g$ denotes Euler's constant.  In all cases, this implies that
$E_{j_1,\snt_1} \le C' j_1\snt_1^{-1}$, for a finite constant~$C'$, implying that the probability 
that the processes $Y\ur$ and~$\hY\ur$ ever differ, if started at time~$\snt_1$ in state~$j_1$, is of order 
$O(j_1\snt_1^{-1})$.  

Now, for the process~$\hY\ur$ started at any time~$\snt_0$ in any state~$j_0$, couple it to a process~$Y\ur$
starting at time~$\snt_1$ in the state~$\hY\ur_{\snt_1}$.  For $2\a > 1$, the probability that the two processes 
ever differ in $\snt \ge \snt_1$ is at most
\[
     \frac{C'}{\snt_1} \ex\{\hY\ur_{\snt_1} \giv \hY\ur_{\snt_0} = j_0\} 
                 \Le \frac{C'}{\snt_1}\,C(j_0 + 2r/(2\a-1))\,\Bigl(\frac{\snt_1}{\snt_0}\Bigr)^{2\a-1},
\]
where, for the last bound, we have used~\Ref{ADB-Yur-mean}, because the mean of~$Y\ur$ grows faster than
that of~$\hY\ur$.  Hence
the probability that the two processes ever differ in $\snt \ge \snt_1$ is at most
\[
         \frac{C'' j_0}{\snt_0^{2\a-1}}\,\frac1{\snt_1^{2(1-\a)}},
\]
and can be chosen to be as small as desired, for fixed initial state~$(j_0,\snt_0)$, by choosing~$\snt_1$ large enough. 
For $2\a < 1$, the comparable bound is of order $O(j_0 \snt_1^{-1})$, and, for $2\a = 1$, it is of order
$O(j_0 \log(\snt_1)\snt_1^{-1})$.
Hence, after some finite random time, the paths of~$\hY\ur$ are exactly those of a process~$Y\ur$.
In consequence, the conclusions of Theorems~\ref{ADB-sXb-recurrence} and~\ref{ADB-recurrence-2} hold for~$\hY\ur$, 
if $\a < (1-q)\log(1/p)$.  If $\a > (1-q)\log(1/p)$, as in Theorem~\ref{aslimit-disc},
$\snt^{-\a} p^{-J_{\snt}} Y_{\snt}\ur$ converges a.s., and, as in Theorem~\ref{CLT-disc},
\[
     \lim_{\snt\to\infty} \pr[\{\log \snt\}^{-1/2}\{\log Y_{\snt}\ur - \bb \log \snt\} \ge y \giv Y_{\snt_0} = x_0] \Eq  
             1 - \Phi(y/v),
\]
where $v^2 = (1-q)\{\log(1/p)\}^2$.

Note that the argument relating $\hY\ur$ to~$Y\ur$ does not greatly involve the detail of the re-wiring
mechanism.  \adbp{For instance, at step $\snt+1$, suppose that each new vertex adds extra connections to a random number of 
other vertices, chosen without replacement from all those other than the vertex being copied, and that double edges
are then merged, as proposed in \cite{bebek2006degree}.  Let the distribution of the number of extra connections 
be denoted by $\hPi^{(\snt)}$, having point 
probabilities~$\hpi^{(\snt)}$ and mean~$r_{\snt}$.  The distribution~$\Pi_k^{(\snt)}$, with point
probabilities~$\pi_k^{(\snt)}$, that represents the number of edges in the copy of a vertex of degree~$k$ at 
step~$\snt + 1$ is not quite the convolution of $\Pi_k$ and~$\hPi^{(\snt)}$, because some of the extra connections
may duplicate edges to neighbours of the copied vertex, but it is easy to check that the difference in total
variation between the two distributions is at most $k r_{\snt}/m$.  The equations 
corresponding to~\Ref{ADB-Q-matrix-rewire}, defining the analogue $Q^{(\snt)}$ of~$Q\urt$, have
$\a$ replaced by $\a_m := \a + (1-\a)r_m/m$, $\hpi\urt$ replaced by~$\hpi^{(\snt)}$ and $\pi_k\urt$ replaced 
by~$\pi_k^{(\snt)}$.
Then, if there is a distribution~$\hPi^{(\infty)}$ with mean~$r^{(\infty)}$ such that 
$\dtv(\hPi^{(\snt)},\hPi^{(\infty)}) = O(m^{-1})$, and if 
$|r_m - r^{(\infty)}| + O(m^{-1})$, the previous arguments can be used to show that the process is asymptotically
equivalent to one associated with the $Q$-matrix~$Q^{(\infty)}$, obtained by formally replacing
$m$ by~$\infty$ in the definitions of~$Q^{(\snt)}$, to which the asymptotic theorems of the paper can 
once more be applied.
}


\end{document}